\documentclass[twocolumn]{autart}    

\usepackage{graphicx}          
\usepackage[tbtags]{amsmath}
\usepackage{color,xcolor}
\usepackage{amssymb}
\usepackage[misc]{ifsym}

\begin{document}

\begin{frontmatter}

\title{Leader-Follower Stochastic Differential Game with Asymmetric Information and Applications\thanksref{footnoteinfo}} 

\thanks[footnoteinfo]{Shi acknowledges the financial support from the National Natural Science Foundations of China (11301011, 11201264, 11571205). Wang acknowledges the financial support from the National Natural Science Foundation for Excellent Young Scholars of China (61422305), the National Natural Science Foundation of China (11371228), the Natural Science Foundation for Distinguished Young Scholars of Shandong Province of China (JQ201418), the Program for New Century Excellent Talents in University of China (NCET-12-0338), the Postdoctoral Foundation of China (2013M540540), and the Research Fund for the Taishan Scholar Project of Shandong Province of China. Xiong acknowledges the financial support from FDCT 076/2012/A3 and MYRG2014-00015-FST. This paper was not presented at any IFAC
meeting. Corresponding author: Guangchen Wang. Tel.: +86 53188392206; Fax: +86 53188365550.
}

\author[JTShi]{Jingtao Shi}\ead{shijingtao@sdu.edu.cn},    
\author[GCWang]{Guangchen Wang}\ead{wguangchen@sdu.edu.cn},               
\author[JXiong]{Jie Xiong}\ead{jiexiong@umac.mo}  

\address[JTShi]{School of Mathematics, Shandong University, Jinan 250100, P.R. China}  
\address[GCWang]{School of Control Science and Engineering, Shandong University, Jinan 250061, P.R. China}             
\address[JXiong]{Department of Mathematics, Faculty of Science and Technology, University of Macau, Macau, P.R. China}        


\begin{abstract}                          
This paper is concerned with a leader-follower stochastic differential game with asymmetric information, where the information available to the follower is based on some sub-$\sigma$-algebra of that available to the leader. Such kind of game problem has wide applications in finance, economics and management engineering such as newsvendor problems, cooperative advertising and pricing problems. Stochastic maximum principles and verification theorems with partial information are obtained, to represent the Stackelberg equilibrium. As applications, a linear-quadratic leader-follower stochastic differential game with asymmetric information is studied. It is shown that the open-loop Stackelberg equilibrium admits a state feedback representation if some system of Riccati equations is solvable.
\end{abstract}

\end{frontmatter}

\section{Introduction}
Throughout this paper, we denote by $\mathbb{R}^n$ the Euclidean space of $n$-dimensional vectors, by $\mathbb{R}^{n\times d}$ the space
of $n\times d$ matrices, by $\mathcal{S}^n$ the space of $n\times n$ symmetric matrices. $\langle\cdot,\cdot\rangle$ and $|\cdot|$ denote the scalar product
and norm in the Euclidean space, respectively. $\top$ appearing in the superscripts denotes the transpose of a matrix. $f_x,f_{xx}$ denote the partial derivative and twice partial derivative with respect to $x$ for a differentiable function $f$.

\subsection{Motivation}

First, we present two examples which motivate us to study leader-follower stochastic differential games with asymmetric information in this paper.

{\it Example 1.1:} (Continuous Time Newsvendor Problem) Let $D(\cdot)$ be the demand rate for some product in a market, which satisfies the stochastic differential equation (SDE)
\begin{equation*}
\left\{
\begin{aligned}
dD(t)=&\ a(\mu-D(t))dt+\sigma dW(t)+\widetilde{\sigma}\widetilde{W}(t),\\
 D(0)=&\ d_0\in\mathbf{R},
\end{aligned}
\right.
\end{equation*}
where $a,\mu,\sigma,\widetilde{\sigma}$ are constants. We consider the market consisting of a manufacturer selling the product to end users through a retailer. At time $t$, the retailer chooses an order rate $q(t)$ for the product and decides its retail price $R(t)$, and is offered a wholesale price $w(t)$ by the manufacturer. We assume that items can be salvaged at unit price $S\geq0$, and that items cannot be stored, i.e., they must be sold instantly or salvaged.

The retailer will obtain an expected profit
\begin{equation*}
\begin{aligned}
 &J_1\big(q(\cdot),R(\cdot),w(\cdot)\big)=\mathbb{E}\int_0^T\big[(R(t)-S)\min[D(t),q(t)]\\
 &\hspace{4.2cm}-(w(t)-S)q(t)\big]dt.
\end{aligned}
\end{equation*}
When the manufacturer has a fixed production cost per unit $M\geq0$, he will get an expected profit
\begin{equation*}
\begin{aligned}
J_2\big(q(\cdot),R(\cdot),w(\cdot)\big)=\mathbb{E}\int_0^T(w(t)-M)q(t)dt.
\end{aligned}
\end{equation*}

Let $\mathcal{F}_t$ denote the $\sigma$-algebra generated by Brownian motions $W(s),\widetilde{W}(s),0\leq s\leq t$. Intuitively $\mathcal{F}_t$ contains all the information up to time $t$. We assume that the information $\mathcal{G}_{1,t},\mathcal{G}_{2,t}$ available to the retailer and the manufacturer at time $t$, respectively, are both sub-$\sigma$-algebras of $\mathcal{F}_t$. Moreover, the information available to them at time $t$ is asymmetric and $\mathcal{G}_{1,t}\subseteq\mathcal{G}_{2,t}$. This can be explained from the practical application's aspect. Specifically, the manufacturer chooses a wholesale price $w(t)$ at time $t$, which is a $\mathcal{G}_{2,t}$-adapted stochastic process. And the retailer chooses an order rate $q(t)$ and a retail price $R(t)$ at time $t$, which are $\mathcal{G}_{1,t}$-adapted stochastic processes. For any $w(\cdot)$, to select a $\mathcal{G}_{1,t}$-adapted process pair $(q^*(\cdot),R^*(\cdot))$ for the retailer such that
\begin{equation*}
\begin{aligned}
       &J_1(q^*(\cdot),R^*(\cdot),w(\cdot))\\
 \equiv&\ J_1\big(q^*(\cdot;w(\cdot)),R^*(\cdot;w(\cdot)),w(\cdot)\big)\\
      =&\ \max\limits_{q(\cdot),R(\cdot)}J_1(q(\cdot),R(\cdot),w(\cdot)),\
\end{aligned}
\end{equation*}
and then to select a $\mathcal{G}_{2,t}$-adapted process $w^*(\cdot)$ for the manufacturer such that
\begin{equation*}
\begin{aligned}
       &J_2(q^*(\cdot),R^*(\cdot),w^*(\cdot))\\
 \equiv&\ J_2\big(q^*(\cdot;w^*(\cdot)),R^*(\cdot;w^*(\cdot)),w^*(\cdot)\big)\\
      =&\ \max\limits_{w(\cdot)}J_2\big(q^*(\cdot;w(\cdot)),R^*(\cdot;w(\cdot)),w(\cdot)\big),
\end{aligned}
\end{equation*}
formulates a leader-follower stochastic differential game with asymmetric information. In this setting, the retailer is the follower and the
manufacturer is the leader. Any process triple $(q^*(\cdot),R^*(\cdot),w^*(\cdot))$ satisfying the above two equalities is called an open-loop Stackelberg equilibrium. In \O ksendal et al. \cite{OSU13}, a time-dependent newsvendor problem with time-delayed information is solved, based on stochastic differential game (with jump-diffusion) approach. But it can not cover our model.

{\it Example 1.2:} (Cooperative Advertising and Pricing Problem) In supply chain management of the market, there are usually two members, the manufacturer and the retailer. Cooperative advertising is an important instrument for aligning manufacturer and retailer decisions. Specifically, we introduce the following SDE, which is the generalization of Sethi's stochastic sales-advertising model introduced by He et al. \cite{HPS09}:
\begin{equation*}
\left\{
\begin{aligned}
dx(t)=&\ \big[\rho u(t)\sqrt{1-x(t)}-\delta x(t)\big]dt\\
      &\ +\sigma(x(t))dW(t)+\widetilde{\sigma}(x(t))d\widetilde{W}(t),\\
 x(0)=&\ x_0\in[0,1],
\end{aligned}
\right.
\end{equation*}
where $x(t)$ represents the awareness share, i.e., the number of aware (or informed) customers expressed as a fraction of the total market at
time $t$, $\rho$ is a response constant, and $\delta$ determines the rate at which potential consumers are lost. $\sigma(x),\widetilde{\sigma}(x)$ are functions satisfying usual conditions.

At time $t\geq0$, the manufacturer decides on the wholesale price $w(t)$ and the cooperative participation rate $\theta(t)$, and the retailer decides the channel's total advertising effort level $u(t)$ and the retail price $p(t)$. The sequence of the events is as follows. At time $t$, first, the manufacturer announces his wholesale price $w(t)$ and participation rate $\theta(t)$. Second, the retailer sets his retail price $p(t)$ and advertising effort rate $u(t)$ as his optimal response to the manufacturer's announced decisions. The retailer accomplishes this by solving an optimization problem to maximize his expected profit
\begin{equation*}
\begin{aligned}
 &\ J_1(w(\cdot),\theta(\cdot),u(\cdot),p(\cdot))\\
=&\ \mathbb{E}\int_0^Te^{-rt}\big[(p(t)-w(t))D(p(t))x(t)\\
 &\qquad\qquad\quad-(1-\theta(t))u^2(t)\big]dt,
\end{aligned}
\end{equation*}
where $r>0$ is the discount rate and $0\leq D(p)\leq 1$ is the demand function satisfying usual conditions. The manufacturer anticipates the retailer's reaction functions and incorporates them into his optimal control problem, and solves for his wholesale price policy $w(t)$ and the participation rate policy $\theta(t)$ at time $t$. Therefore, the manufacturer's optimization problem is to maximize his expected profit
\begin{equation*}
\begin{aligned}
 &\ J_2(w(\cdot),\theta(\cdot),u(\cdot),p(\cdot))\\
=&\ \mathbb{E}\int_0^Te^{-rt}\big[(w(t)-c)D(p(t))x(t)-\theta(t)u^2(t)\big]dt,
\end{aligned}
\end{equation*}
where $c\geq0$ is the constant unit production cost.

Define $\mathcal{F}_t:=\sigma\{W(s),0\leq s\leq t\}\bigvee\sigma\{\widetilde{W}(s),0\leq s\leq t\}$. In the game setting, we assume that the information $\mathcal{G}_{1,t},\mathcal{G}_{2,t}$ available to the retailer and the manufacturer at time $t$, respectively, are both sub-$\sigma$-algebras of the complete information filtration $\mathcal{F}_t$. Moreover, the information available at time $t$ to them is asymmetric and $\mathcal{G}_{1,t}\subseteq\mathcal{G}_{2,t}$. In detail, the wholesale price $w(\cdot)$ and the participation rate $\theta(\cdot)$ of the manufacturer are $\mathcal{G}_{2,t}$-adapted processes. For the retailer, his advertising effort level $u(\cdot)$ and retail price $p(\cdot)$ are to be $\mathcal{G}_{1,t}$-adapted processes. For any $(w(\cdot),\theta(\cdot))$, first, to select a suitable process pair $(u^*(\cdot),p^*(\cdot))$ for the retailer such that
\begin{equation*}
\begin{aligned}
      &\ J_1(w(\cdot),\theta(\cdot),u^*(\cdot),p^*(\cdot))\\
\equiv&\ J_1\big(w(\cdot),\theta(\cdot),u^*(\cdot;(w(\cdot),\theta(\cdot)),p^*(\cdot;(w(\cdot),\theta(\cdot))\big)\\
     =&\ \max\limits_{u(\cdot),p(\cdot)\geq0}J_1(w(\cdot),\theta(\cdot),u(\cdot),p(\cdot)),
\end{aligned}
\end{equation*}
and then to select a suitable process pair $(w^*(\cdot),\theta^*(\cdot))$ for the manufacturer such that
{\small\begin{equation*}
\begin{aligned}
      &\ J_2(w^*(\cdot),\theta^*(\cdot),u^*(\cdot),p^*(\cdot))\\
\equiv&\ J_2\big(w^*(\cdot),\theta^*(\cdot),u^*(\cdot;(w^*(\cdot),\theta^*(\cdot)),p^*(\cdot;(w^*(\cdot),\theta^*(\cdot))\big)\\
     =&\ \max\limits_{w(\cdot),0\leq\theta(\cdot)\leq1}J_2\big(w(\cdot),\theta(\cdot),u^*(\cdot;(w(\cdot),\theta(\cdot)),\\
      &\qquad\qquad\qquad\quad p^*(\cdot;(w(\cdot),\theta(\cdot))\big),
\end{aligned}
\end{equation*}}
formulates a leader-follower stochastic differential game with asymmetric information. See also \cite{HPS09} for more details about cooperative advertising and pricing models in dynamic stochastic supply chains, where feedback Stackelberg equilibrium is obtained applying dynamic programming approach for stochastic differential game. However, the asymmetric information was not considered there.

\subsection{Problem Formulation}

Inspired by the examples above, we study leader-follower stochastic differential games with asymmetric information in this paper.

Let $0<T<\infty$ be a finite time duration and $(\Omega,\mathcal{F},\mathbb{P})$ be a complete probability space.
$(W(\cdot),\widetilde{W}(\cdot))$ is a standard $\mathbb{R}^{d_1+d_2}$-valued Brownian motion. Let $\{\mathcal{F}_t\}_{0\leq t\leq T}$ be the natural augmented filtration generated by $(W(\cdot),\widetilde{W}(\cdot))$ and $\mathcal{F}_T=\mathcal{F}$.

Suppose that the state of the system is described by the SDE
\begin{equation}\label{state equation for the follower}
\left\{
\begin{aligned}
dx^{u_1,u_2}(t)&=b\big(t,x^{u_1,u_2}(t),u_1(t),u_2(t)\big)dt\\
               &\quad+\sigma\big(t,x^{u_1,u_2}(t),u_1(t),u_2(t)\big)dW(t)\\
               &\quad+\widetilde{\sigma}\big(t,x^{u_1,u_2}(t),u_1(t),u_2(t)\big)d\widetilde{W}(t),\\
 x^{u_1,u_2}(0)&=\ x_0,
\end{aligned}
\right.
\end{equation}
where $u_1(\cdot)$ and $u_2(\cdot)$ are control processes taken by the two players in the game, labeled 1 (the follower) and 2 (the leader), with values in nonempty convex sets $U_1\subseteq\mathbb R^{m_1}$, $U_2\subseteq\mathbb R^{m_2}$, respectively. $x^{u_1,u_2}(\cdot)$, the solution to SDE (\ref{state equation for the follower}) with values in $\mathbb{R}^n$, is the corresponding state process with initial state $x_0\in\mathbb{R}^n$. Here $b(t,x,u_1,u_2):\Omega\times[0,T]\times\mathbb{R}^n\times U_1\times U_2\rightarrow\mathbb{R}^n,\ \sigma(t,x,u_1,u_2):\Omega\times[0,T]\times\mathbb{R}^n\times U_1\times U_2\rightarrow\mathbb{R}^{n\times d_1},\ \widetilde{\sigma}(t,x,u_1,u_2):\Omega\times[0,T]\times\mathbb{R}^n\times U_1\times U_2\rightarrow\mathbb{R}^{n\times d_2}$ are given $\mathcal{F}_t$-adapted processes, for each $(x,u_1,u_2)$.

Let us now explain the asymmetric information feature between the follower (player 1) and the leader (player 2) in this paper. Player 1 is the follower, which means that the information available to him at time $t$ is based on some sub-$\sigma$-algebra $\mathcal{G}_{1,t}\subseteq\mathcal{G}_{2,t}$, where $\mathcal{G}_{2,t}$ is the information available to the leader at time $t$. We assume in this and next sections that both $\mathcal{G}_{1,t}$ and $\mathcal{G}_{2,t}$ are sub-$\sigma$-algebras of the complete information filtration $\mathcal{F}_t$. That is, we have $\mathcal{G}_{1,t}\subseteq\mathcal{G}_{2,t}\subseteq\mathcal{F}_t$. We define the admissible control sets of the follower and the leader, respectively, as follows.
\begin{equation}\label{admissible control set for the follower}
\begin{aligned}
\mathcal{U}_1:=&\ \Big\{u_1\big|u_1:\Omega\times[0,T]\rightarrow U_1\mbox{ is }\mathcal{G}_{1,t}\mbox{-adapted}\\
&\quad\mbox{and }\sup\limits_{0\leq t\leq T}\mathbb{E}|u_1(t)|^i<\infty,\ i=1,2,\cdots\Big\},
\end{aligned}
\end{equation}
\begin{equation}\label{admissible control set for the leader}
\begin{aligned}
\mathcal{U}_2:=&\ \Big\{u_2\big|u_2:\Omega\times[0,T]\rightarrow U_2\mbox{ is }\mathcal{G}_{2,t}\mbox{-adapted}\\
&\quad\mbox{and }\sup\limits_{0\leq t\leq T}\mathbb{E}|u_2(t)|^i<\infty,\ i=1,2,\cdots\Big\}.
\end{aligned}
\end{equation}

In the game problem, knowing that the leader has chosen $u_2(\cdot)\in\mathcal{U}_2$, the follower would like to choose a $\mathcal{G}_{1,t}$-adapted control $u_1^*(\cdot)=u_1^*(\cdot;u_2(\cdot))$ to minimize his cost functional
\begin{equation}\label{cost functional for the follower}
\begin{aligned}
 &\ J_1(u_1(\cdot),u_2(\cdot))\\
=&\ \mathbb{E}\bigg[\int_0^Tg_1\big(t,x^{u_1,u_2}(t),u_1(t),u_2(t)\big)dt\\
 &\quad+G_1(x^{u_1,u_2}(T))\bigg].
\end{aligned}
\end{equation}
Here $g_1(t,x,u_1,u_2):\Omega\times[0,T]\times\mathbb{R}^n\times U_1\times U_2\rightarrow\mathbb{R}$ is an $\mathcal{F}_t$-adapted process, and $G_1(x):\Omega\times\mathbb{R}^n\rightarrow\mathbb{R}$ is an $\mathcal{F}_T$-measurable random variable, for each $(x,u_1,u_2)$. Now the follower encounters a stochastic optimal control problem with partial information.

{\bf Problem of the follower.}\ For any chosen $u_2(\cdot)\in\mathcal{U}_2$ by the leader, choose a $\mathcal{G}_{1,t}$-adapted control $u_1^*(\cdot)=u_1^*(\cdot;u_2(\cdot))\in\mathcal{U}_1$, such that
\begin{equation}\label{problem of the follower}
\begin{aligned}
J_1(u_1^*(\cdot),u_2(\cdot))\equiv&\ J_1\big(u_1^*(\cdot;u_2(\cdot)),u_2(\cdot)\big)\\
                                 =&\ \inf\limits_{u_1\in\mathcal{U}_1}J_1(u_1(\cdot),u_2(\cdot)),
\end{aligned}
\end{equation}
subject to (\ref{state equation for the follower}) and (\ref{cost functional for the follower}). Such a $u_1^*(\cdot)=u_1^*(\cdot;u_2(\cdot))$ is called a (partial information) optimal control, and the corresponding solution $x^{u_1^*,u_2}(\cdot)$ to (\ref{state equation for the follower}) is called a (partial information) optimal state process for the follower.

In the following procedure of the game problem, once knowing that the follower would take such an optimal control $u_1^*(\cdot)=u_1^*(\cdot;u_2(\cdot))$, the leader would like to choose a $\mathcal{G}_{2,t}$-adapted control $u_2^*(\cdot)$ to minimize his cost functional
\begin{equation}\label{cost functional for the leader-knowing the optimal control of the follower}
\begin{aligned}
 &\ J_2(u_1^*(\cdot),u_2(\cdot))\\
=&\ \mathbb{E}\bigg[\int_0^Tg_2\big(t,x^{u_1^*,u_2}(t),u_1^*(t;u_2(t)),u_2(t)\big)dt\\
 &\quad+G_2(x^{u_1^*,u_2}(T))\bigg].
\end{aligned}
\end{equation}
Here $g_2(t,x,u_1,u_2):\Omega\times[0,T]\times\mathbb{R}^n\times U_1\times U_2\rightarrow\mathbb{R},G_2(x):\Omega\times\mathbb{R}^n\rightarrow\mathbb{R}$ are given $\mathcal{F}_t$-adapted processes, for each $(x,u_1,u_2)$. Now the leader encounters a stochastic optimal control problem with partial information.

{\bf Problem of the leader.}\ Find a $\mathcal{G}_{2,t}$-adapted control $u_2^*(\cdot)\in\mathcal{U}_2$, such that
\begin{equation}\label{problem of the leader}
\begin{aligned}
J_2(u_1^*(\cdot),u_2^*(\cdot))=&\ J_2\big(u_1^*(\cdot;u_2^*(\cdot)),u_2^*(\cdot)\big)\\
                              =&\ \inf\limits_{u_2\in\mathcal{U}_2}J_2\big(u_1^*(\cdot;u_2(\cdot)),u_2(\cdot)\big),
\end{aligned}
\end{equation}
subject to (\ref{state equation for the follower}) and (\ref{cost functional for the leader-knowing the optimal control of the follower}). Such a $u_2^*(\cdot)$ is called a (partial information) optimal control, and the corresponding solution $x^*(\cdot)\equiv x^{u_1^*,u_2^*}(\cdot)$ to (\ref{state equation for the follower}) is called a (partial information) optimal state process for the leader. We will restate the problem for the leader in more detail in the next section, since its precise description has to involve the solution to the {\bf Problem of the follower}.

We refer to the problem mentioned above as a {\it leader-follower stochastic differential game with asymmetric information}. If there exists a control process pair $(u_1^*(\cdot),u_2^*(\cdot))=\big(u_1^*(\cdot;u_2^*(\cdot)),u_2^*(\cdot)\big)$ satisfying (\ref{problem of the follower}) and (\ref{problem of the leader}), we refer to it as an {\it open-loop Stackelberg equilibrium}.

In this paper, we impose the following assumptions.

\noindent{\bf (A1.1)}\ {\it For each $\omega\in\Omega$, the functions $b,\sigma,\widetilde{\sigma},g_1$ are twice continuously differentiable with respect to $(x,u_1,u_2)$. For each $\omega\in\Omega$, functions $g_2$ and $G_1,G_2$ are continuously differentiable with respect to $(x,u_1,u_2)$ and $x$, respectively. Moreover, for each $\omega\in\Omega$ and any $(t,x,u_1,u_2)\in[0,T]\times\mathbb{R}^n\times\mathbb{R}^{m_1}\times\mathbb{R}^{m_2}$, there exists $C>0$ such that
\begin{equation*}
\begin{aligned}
&\big(1+|x|+|u_1|+|u_2|\big)^{-1}\big|\phi(t,x,u_1,u_2)\big|\\
&+\big|\phi_x(t,x,u_1,u_2)\big|+\big|\phi_{u_1}(t,x,u_1,u_2)\big|\\
&+\big|\phi_{u_2}(t,x,u_1,u_2)\big|+\big|\phi_{xx}(t,x,u_1,u_2)\big|\\
&+\big|\phi_{u_1u_1}(t,x,u_1,u_2)\big|+\big|\phi_{u_2u_2}(t,x,u_1,u_2)\big|\leq C,
\end{aligned}
\end{equation*}
for $\phi=b,\sigma,\widetilde{\sigma}$, and
\begin{equation*}
\begin{aligned}
&\big(1+|x|^2\big)^{-1}\big|G_1(x)\big|+\big(1+|x|\big)^{-1}\big|G_{1x}(x)\big|\\
&+\big(1+|x|^2\big)^{-1}\big|G_2(x)\big|+\big(1+|x|\big)^{-1}\big|G_{2x}(x)\big|\leq C,
\end{aligned}
\end{equation*}
\begin{equation*}
\begin{aligned}
&\big(1+|x|^2+|u_1|^2+|u_2|^2\big)^{-1}\big|g_1(t,x,u_1,u_2)\big|\\
&+\big(1+|x|+|u_1|+|u_2|\big)^{-1}\Big(\big|g_{1x}(t,x,u_1,u_2)\big|\\
&+\big|g_{1u_1}(t,x,u_1,u_2)\big|+\big|g_{1u_2}(t,x,u_1,u_2)\big|\Big)\\
&+\big|g_{1xx}(t,x,u_1,u_2)\big|+\big|g_{1u_1u_1}(t,x,u_1,u_2)\big|\\
&+\big|g_{1u_2u_2}(t,x,u_1,u_2)\big|\leq C,
\end{aligned}
\end{equation*}
\begin{equation*}
\begin{aligned}
&\big(1+|x|^2+|u_1|^2+|u_2|^2\big)^{-1}\big|g_2(t,x,u_1,u_2)\big|\\
&+\big(1+|x|+|u_1|+|u_2|\big)^{-1}\Big(\big|g_{2x}(t,x,u_1,u_2)\big|\\
&+\big|g_{2u_1}(t,x,u_1,u_2)\big|+\big|g_{2u_2}(t,x,u_1,u_2)\big|\Big)\leq C.
\end{aligned}
\end{equation*}}

\subsection{Literature Review and Contributions of This Paper}

Initiated by Issacs \cite{Isaacs54}, differential games are useful in modeling dynamic systems where more than one decision maker are involved. Differential games have been investigated by many authors and have been found to be a useful tool in many applications, particularly in biology, economics and finance. Stochastic differential games are differential games for stochastic systems involving noise terms. See the monographs by Basar and Olsder \cite{BO82} for more information about differential games. For some most recent developments for stochastic differential games and their applications, please refer to Yong \cite{Yong90}, Hamad\`{e}ne \cite{Ha99}, Wu \cite{Wu05}, An and \O ksendal \cite{AO08}, Buckdahn and Li \cite{BL08}, Wang and Yu \cite{WY10,WY12}, Yu \cite{Yu12}, Hui and Xiao \cite{HX12,HX14}, Shi \cite{Shi13} and the references therein.

Leader-follower stochastic differential game is the dynamic and stochastic formulation of the well-known Stackelberg game, which was introduced by Stackelberg \cite{Stackelberg34} in 1934, when he defined a concept of a hierarchical solution for markets where some firms have power of domination over others. This solution concept is now known as the Stackelberg equilibrium which, in the context of two-person nonzero-sum games, involves players with asymmetric roles, one leader and one follower. Early study for stochastic Stackelberg differential games can be seen in Basar \cite{Basar79}. In detail, a leader-follower stochastic differential game (or stochastic Stackelberg differential game) proceeds with the follower aims at minimizing his cost functional in accordance with the leader's strategy on the whole duration of the game. Anticipating the follower's optimal response depending on his entire strategy, the leader chooses an optimal one in advance to minimize his own cost functional, based on the stochastic Hamiltonian system satisfied by the follower's optimal response. The pair of the leader's optimal strategy and the follower's optimal response is known as the Stackelberg equilibrium.

To our best knowledge, there are few papers on leader-follower stochastic differential games. A pioneer work was done by Yong \cite{Yong02}, where a linear-quadratic (LQ) leader-follower stochastic differential game was introduced and studied. The coefficients of the system and the cost functionals are random, the controls enter the diffusion term of the state equation, and the weight matrices for the controls in the cost functionals are not necessarily positive definite. To give a state feedback representation of the open-loop Stackelberg equilibrium (in a non-anticipating way), the related Riccati equations are derived and sufficient conditions for the existence of their solution with deterministic coefficients are discussed. Bensoussan et al. \cite{BCS12} obtained the global maximum principles for both open-loop and closed-loop stochastic Stackelberg differential games, whereas the diffusion term does not contain the controls. The solvability of related Riccati equations is discussed, in order to obtain the state feedback Stackelberg equilibrium.

In this paper, we initiate to study a leader-follower stochastic differential game with asymmetric information, which distinguishes itself from the literatures mentioned above in the following aspects. (i) In our framework, both information filtration available to the leader and the follower could be sub-$\sigma$-algebras of the complete information filtration naturally generated by the random noise source. Moreover, the information available to the follower is based on some sub-$\sigma$-algebra of that available to the leader. This gives a new explanation for the asymmetric information feature between the follower and the leader, and endows our problem formulation more practical meanings in reality. (ii) An important class of LQ leader-follower stochastic differential games with asymmetric information is first proposed and then solved. It consists of a stochastic optimal control problem of SDE with partial observation for the follower, and followed by a stochastic optimal control problem of conditional mean-field forward-backward stochastic differential equation (FBSDE) with complete information for the leader. This problem is new in differential game theory and have considerable impacts in both theoretical analysis and practical meaning with future application prospect, although it has intrinsic mathematical difficulties. Note that in \cite{Yong02}, both problems for the follower and leader are to be solved stochastic optimal control problems with complete information. (iii) The open-loop Stackelberg equilibrium of this LQ game problem with asymmetric information, is characterized in terms of the forward-backward stochastic differential filtering equation (FBSDFE), which arises naturally in our setup. To our best knowledge, these FBSDFEs are new in both stochastic control and filtering theory. In particular, they are different from those in Zhang \cite{Zhang90}, Huang et al. \cite{HWX09} and Wang and Yu \cite{WY12}. (iv) State feedback representations for the optimal controls of the follower and the leader, are explicitly given with the help of some new Riccati equations. Both problems for the follower and the leader have mathematical difficulties, and we overcome them via some measure transformation and filtering technique in Xiong \cite{Xiong08}, Wang et al. \cite{WWX13}, linear FBSDE decoupling technique in \cite{Yong02} and mean-field FBSDE decoupling technique in Yong \cite{Yong13}, respectively.

The rest of this paper is organized as follows. In Section 2, we formulate the problems of the follower and the leader, and then solve them in this order. In Subsection 2.1, we first solve the stochastic optimal control problem with partial information of the follower, applying the theory of controlled SDEs with partial information. The partial information maximum principle and verification theorem are shown. In Subsection 2.2, we then formulate the stochastic optimal control problem of conditional mean-field FBSDE with partial information of the leader, regarding the stochastic Hamiltonian system corresponding to the follower's optimal response as his state equation. Via controlled conditional mean-field FBSDEs, the maximum principle and verification theorem with partial information are both given. In Section 3, we apply our theoretical results to an LQ leader-follower stochastic differential game with asymmetric information. Specifically, Subsection 3.1 is devoted to the solution of an LQ stochastic optimal control problem with partial observation of the follower. The key technique is to get some observable optimal controls by explicitly computing the optimal filtering estimates of the corresponding adjoint processes, when applying the maximum principle approach and measure transformation. The state feedback form of the observable optimal control for the follower is then represented by some symmetric Riccati equation under appropriate assumptions. Here, the state of the follower is represented by some FBSDFE. Subsection 3.2 is devoted to the solution of an LQ stochastic optimal control problem of conditional mean-field FBSDE with complete information of the leader. The state feedback representation of optimal control is derived. Since the filtering estimates of the original state and the control of the leader are involved in its state equation, we encounter difficulty. We overcome this difficulty by applying the maximum principle and verification theorem for conditional mean-field FBSDEs, via the solution to some adjoint conditional mean-field FBSDE. The state feedback representation for the optimal control of the leader (together with the optimal control of the follower) is obtained via some new system of Riccati equations in the double dimensional spaces, when the dimensional augmentation by Yong \cite{Yong02} is applied. Thus, the state feedback representation of the open-loop Stackelberg equilibrium can be given. Finally, Section 4 gives some concluding remarks.

\section{Maximum Principle and Verification Theorem for Stackelberg Equilibrium}

\subsection{The Follower's Problem}

In this subsection, we first solve the partial information stochastic optimal control problem of the follower, that is {\bf Problem of the follower}.

For any chosen $u_2(\cdot)\in\mathcal{U}_2$, suppose that there exists an optimal control $u_1^*(\cdot)$ for the follower, and the corresponding optimal state $x^{u_1^*,u_2}(\cdot)$ is the solution to (\ref{state equation for the follower}). Let an $\mathcal{F}_t$-adapted process triple $(q(\cdot),k(\cdot),\widetilde{k}(\cdot))\in\mathbb{R}^n\times\mathbb{R}^{n\times d_1}\times\mathbb{R}^{n\times d_2}$ uniquely solves the adjoint BSDE of the follower
\begin{equation}\label{adjoint equation for the follower}
\left\{
\begin{aligned}
-dq(t)=&\ \Big\{b_x\big(t,x^{u_1^*,u_2}(t),u_1^*(t),u_2(t)\big)q(t)\\
       &\ +\sum\limits_{j=1}^{d_1}\sigma_x^j\big(t,x^{u_1^*,u_2}(t),u_1^*(t),u_2(t)\big)k^j(t)\\
       &\ +\sum\limits_{j=1}^{d_2}\widetilde{\sigma}_x^j\big(t,x^{u_1^*,u_2}(t),u_1^*(t),u_2(t)\big)\widetilde{k}^j(t)\\
       &\ -g_{1x}\big(t,x^{u_1^*,u_2}(t),u_1^*(t),u_2(t)\big)\Big\}dt\\
       &-k(t)dW(t)-\widetilde{k}(t)d\widetilde{W}(t),\\
  q(T)=&-G_{1x}(x^{u_1^*,u_2}(T)),
\end{aligned}
\right.
\end{equation}
where we denote $\sigma_x\equiv\big(\sigma_x^1,\sigma_x^2,\cdots,\sigma_x^{d_1}\big)^\top, \widetilde{\sigma}_x\equiv\big(\widetilde{\sigma}_x^1,\widetilde{\sigma}_x^2,\cdots,\widetilde{\sigma}_x^{d_2}\big)^\top$, and similar notations will be used.

Define the Hamiltonian function of the follower $H_1:\Omega\times[0,T]\times\mathbb{R}^n\times U_1\times U_2\times\mathbb{R}^n\times\mathbb{R}^{n\times d_1}\times\mathbb{R}^{n\times d_2}\rightarrow\mathbb{R}$ as
\begin{equation}\label{Hamiltonian function for the follower}
\begin{aligned}
 &H_1\big(t,x,u_1,u_2;q,k,\widetilde{k}\big)=\big\langle q,b(t,x,u_1,u_2)\big\rangle\\
 &\ +\mbox{tr}\big\{k^\top\sigma(t,x,u_1,u_2)\big\}+\mbox{tr}\big\{\widetilde{k}^\top\widetilde{\sigma}(t,x,u_1,u_2)\big\}\\
 &\ -g_1(t,x,u_1,u_2).
\end{aligned}
\end{equation}

Similar to \cite{WY12}, we have the following results.

\noindent{\bf Proposition 2.1 (Partial Information Maximum Principle)}
\quad {\it Suppose that {\bf (A1.1)} holds. For any given $u_2(\cdot)\in\mathcal{U}_2$, let $u_1^*(\cdot)$ be the optimal control for {\bf Problem of the follower}, and $x^{u_1^*,u_2}(\cdot)$ be the corresponding optimal state. Let $(q(\cdot),k(\cdot),\widetilde{k}(\cdot))$ be the adjoint process triple. Then we have
\begin{equation}\label{maximum condition for the follower-necessary MP}
\begin{aligned}
&\mathbb{E}\Big[\big\langle H_{1u_1}\big(t,x^{u_1^*,u_2}(t),u_1^*(t),u_2(t);q(t),k(t),\widetilde{k}(t)\big),\\
&\qquad u_1-u_1^*(t)\big\rangle\Big|\mathcal{G}_{1,t}\Big]\geq0,\ a.e.t\in[0,T], a.s.,
\end{aligned}
\end{equation}
holds, for any $u_1\in U_1$.}

\vspace{1mm}

\noindent{\bf Proposition 2.2 (Partial Information Verification Theorem)}
\quad {\it Suppose that {\bf (A1.1)} holds. For any given $u_2(\cdot)$, let $u_1^*(\cdot)\in\mathcal{U}_1$ and $x^{u_1^*,u_2}(\cdot)$ be the corresponding state. Let $(q(\cdot),k(t),\widetilde{k}(\cdot))$ be the adjoint process triple. Moreover, for each $(t,\omega)\in[0,T]\times\Omega$, $H_1\big(t,\cdot,\cdot,u_2(t);q(t),k(t),\widetilde{k}(t)\big)$
is concave, $G_1(\cdot)$ is convex, and
{\small\begin{equation}\label{maximum condition for the follower-sufficient MP}
\begin{aligned}
&\mathbb{E}\Big[H_1\big(t,x^{u_1^*,u_2}(t),u_1^*(t),u_2(t);q(t),k(t),\widetilde{k}(t)\big)\Big|\mathcal{G}_{1,t}\Big]\\
=&\max\limits_{u_1\in U_1}\mathbb{E}\Big[H_1\big(t,x^{u_1^*,u_2}(t),u_1,u_2(t);q(t),k(t),\widetilde{k}(t)\big)\Big|\mathcal{G}_{1,t}\Big],
\end{aligned}
\end{equation}}
holds for $a.e.t\in[0,T]$, a.s. Then $u^*_1(\cdot)$ is an optimal control for {\bf Problem of the follower}.}

\subsection{The Leader's Problem}

In this subsection, we first state the partial information stochastic optimal control problem of the leader in detail, that is {\bf Problem of the leader}. Then we give the maximum principle and verification theorem for it.

For any $u_2(\cdot)\in\mathcal{U}_2$, by the maximum condition (\ref{maximum condition for the follower-necessary MP}), we assume that a functional $u_1^*(t)=u_1^*\big(t;\hat{x}^{u_1^*,\hat{u}_2}(t),\hat{u}_2(t),\\\hat{q}(t),\hat{k}(t),\hat{\widetilde{k}}(t)\big)$ is uniquely defined, where
\begin{equation}\label{filtering estimates}
\left\{
\begin{aligned}
   &\hat{x}^{u_1^*,\hat{u}_2}(t):=\mathbb{E}\big[x^{u_1^*,u_2}(t)\big|\mathcal{G}_{1,t}\big],\\
   &\hat{u}_2(t):=\mathbb{E}\big[u_2(t)\big|\mathcal{G}_{1,t}\big],\quad\hat{q}(t):=\mathbb{E}\big[q(t)\big|\mathcal{G}_{1,t}\big],\\
   &\hat{k}(t):=\mathbb{E}\big[k(t)\big|\mathcal{G}_{1,t}\big],\quad\hat{\widetilde{k}}(t):=\mathbb{E}\big[\widetilde{k}(t)\big|\mathcal{G}_{1,t}\big].
\end{aligned}
\right.
\end{equation}
For the simplicity of notations, we denote $x^{u_2}(\cdot)\equiv x^{u_1^*,u_2}(\cdot)$ and define $\phi^L$ on $\Omega\times[0,T]\times\mathbb{R}^n\times U_2$ as
\begin{equation*}
\begin{aligned}
&\phi^L\big(t,x^{u_2}(t),u_2(t)\big):=\phi\big(t,x^{u_1^*,u_2}(t),\\
&\qquad u_1^*\big(t;\hat{x}^{u_1^*,\hat{u}_2}(t),\hat{u}_2(t),\hat{q}(t),\hat{k}(t),\hat{\widetilde{k}}(t)\big),u_2(t)\big),
\end{aligned}
\end{equation*}
for $\phi=b,\sigma,\widetilde{\sigma},g_1$, respectively. Then after substituting the above control process $u_1^*(\cdot)$ into the follower's adjoint equation (\ref{adjoint equation for the follower}), the leader encounters the controlled FBSDE system
\begin{equation}\label{state equation for the leader}
\left\{
\begin{aligned}
dx^{u_2}(t)=&\ b^L\big(t,x^{u_2}(t),u_2(t)\big)dt\\
            &+\sigma^L\big(t,x^{u_2}(t),u_2(t)\big)dW(t)\\
            &+\widetilde{\sigma}^L\big(t,x^{u_2}(t),u_2(t)\big)d\widetilde{W}(t),\\
     -dq(t)=&\ \Big\{b_x^L\big(t,x^{u_2}(t),u_2(t)\big)q(t)\\
            &\ +\sum\limits_{j=1}^{d_1}\sigma_x^{Lj}\big(t,x^{u_2}(t),u_2(t)\big)k^j(t)\\
            &\ +\sum\limits_{j=1}^{d_2}\widetilde{\sigma}_x^{Lj}\big(t,x^{u_2}(t),u_2(t)\big)\widetilde{k}^j(t)\\
            &\ -g_{1x}^L\big(t,x^{u_2}(t),u_2(t)\big)\Big\}dt\\
            &-k(t)dW(t)-\widetilde{k}(t)d\widetilde{W}(t),\\
 x^{u_2}(0)=&\ x_0,\ q(T)=-G_{1x}(x^{u_2}(T)).
\end{aligned}
\right.
\end{equation}
Note that (\ref{state equation for the leader}) is a controlled {\it conditional mean-field FBSDE}, which now is regarded as the ``state" equation of the leader. That is to say, the state for the leader is the quadruple $(x^{u_2}(\cdot),q(\cdot),k(\cdot),\widetilde{k}(\cdot))$.

\noindent{\bf Remark 2.2}\quad The equality $u_1^*(t)=u_1^*\big(t;\hat{x}^{u_1^*,\hat{u}_2}(t),\hat{u}_2(t),\\\hat{q}(t),\hat{k}(t),\hat{\widetilde{k}}(t)\big)$ does not hold in general. However, for LQ case, it is satisfied and we will make this point clear in the next section.

Define
\begin{equation}\label{cost functional for the leader}
\begin{aligned}
&J_2^L(u_2(\cdot)):=J_2(u_1^*(\cdot),u_2(\cdot))\\
&=\mathbb{E}\bigg[\int_0^Tg_2\big(t,x^{u_1^*,u_2}(t),u_1^*(t),u_2(t)\big)dt\\
&\qquad+G_2(x^{u_1^*,u_2}(T))\bigg]\\
&\equiv\mathbb{E}\bigg[\int_0^Tg_2\big(t,x^{u_1^*,u_2}(t),u_1^*\big(t;\hat{x}^{u_1^*,\hat{u}_2}(t),\hat{u}_2(t),\\
&\qquad\ \hat{q}(t),\hat{k}(t),\hat{\widetilde{k}}(t)\big),u_2(t)\big)dt+G_2(x^{u_1^*,u_2}(T))\bigg]\\
&:=\mathbb{E}\bigg[\int_0^Tg_2^L\big(t,x^{u_2}(t),u_2(t)\big)dt+G_2(x^{u_2}(T))\bigg],
\end{aligned}
\end{equation}
where $g_2^L:\Omega\times[0,T]\times\mathbb{R}^n\times U_2\rightarrow\mathbb{R}$. Note the cost functional of the leader is also conditional mean-field's type. We propose the stochastic optimal control problem with partial information of the leader as follows.

{\bf Problem of the leader.}\ Find a $\mathcal{G}_{2,t}$-adapted control $u_2^*(\cdot)\in\mathcal{U}_2$, such that
\begin{equation}\label{problem of the leader-another form}
J_2^L(u_2^*(\cdot))=\inf\limits_{u_2\in\mathcal{U}_2}J_2^L(u_2(\cdot)),
\end{equation}
subject to (\ref{state equation for the leader}) and (\ref{cost functional for the leader}). Such a $u_2^*(\cdot)$ is called a (partial information) optimal control, and the corresponding solution $x^*(\cdot)\equiv x^{u_2^*}(\cdot)$ to (\ref{state equation for the leader}) is called a (partial information) optimal state process for the leader.

Suppose that there exists an optimal control $u_2^*(\cdot)$ for the leader, and the corresponding state
$(x^*(\cdot),q^*(\cdot),k^*(\cdot),\\\widetilde{k}^*(\cdot))$ is the solution to (\ref{state equation for the leader}). Define the Hamiltonian function of the leader $H_2:\Omega\times[0,T]\times\mathbb{R}^n\times U_2\times\mathbb{R}^n\times\mathbb{R}^{n\times d_1}\times\mathbb{R}^{n\times d_2}\times\mathbb{R}^n\times\mathbb{R}^n\times\mathbb{R}^{n\times d_1}\times\mathbb{R}^{n\times d_2}\times\mathbb{R}^n\rightarrow\mathbb{R}$ as
\begin{equation}\label{Hamiltonian function for the leader}
\begin{aligned}
&H_2\big(t,x^{u_2},u_2,q,k,\widetilde{k};y,z,\widetilde{z},p\big)\\
&=\big\langle y,b^L(t,x^{u_2},u_2)\big\rangle+\mbox{tr}\big\{z^\top\sigma^L(t,x^{u_2},u_2)\big\}\\
&\quad+\mbox{tr}\big\{\widetilde{z}^\top\widetilde{\sigma}^L(t,x^{u_2},u_2)\big\}+g_2^L(t,x^{u_2},u_2)\\
&\quad-\Big\langle p,b_x^L(t,x^{u_2},u_2)q+\sum\limits_{j=1}^{d_1}\sigma_x^{Lj}(t,x^{u_2},u_2)k^j\\
&\qquad\quad+\sum\limits_{j=1}^{d_2}\widetilde{\sigma}_x^{Lj}(t,x^{u_2},u_2)\widetilde{k}^j-g_{1x}^L(t,x^{u_2},u_2)\Big\rangle.
\end{aligned}
\end{equation}
Let $(y(\cdot),z(\cdot),\widetilde{z}(\cdot),p(\cdot))\in\mathbb{R}^n\times\mathbb{R}^{n\times d_1}\times\mathbb{R}^{n\times d_2}\times\mathbb{R}^n$ be the unique $\mathcal{F}_t$-adapted solution to the adjoint conditional mean-field FBSDE of the leader
\begin{equation}\label{adjoint equation for the leader}
\left\{
\begin{aligned}
 dp(t)=&\ \Big\{b_x^{L*}(t)p(t)+\mathbb{E}\big[b_{\hat{x}}^{L*}(t)p(t)\big|\mathcal{G}_{1,t}\big]\Big\}dt\\
       &\ +\sum\limits_{j=1}^{d_1}\Big\{\sigma_x^{L*j}(t)p(t)\\
       &\qquad\quad+\mathbb{E}\big[\sigma_{\hat{x}}^{L*j}(t)p(t)\big|\mathcal{G}_{1,t}\big]\Big\}dW^j(t)\\
       &\ +\sum\limits_{j=1}^{d_2}\Big\{\widetilde{\sigma}_x^{L*j}(t)p(t)\\
       &\qquad\quad+\mathbb{E}\big[\widetilde{\sigma}_{\hat{x}}^{L*j}(t)p(t)\big|\mathcal{G}_{1,t}\big]\Big\}d\widetilde{W}^j(t),\\
-dy(t)=&\ \bigg\{b_x^{L*}(t)y(t)+\mathbb{E}\big[b_{\hat{x}}^{L*}(t)y(t)\big|\mathcal{G}_{1,t}\big]\\
       &+\sum\limits_{j=1}^{d_1}\Big[\sigma_x^{L*j}(t)z^j(t)+\mathbb{E}\big[\sigma_{\hat{x}}^{L*j}(t)z^j(t)\big|\mathcal{G}_{1,t}\big]\Big]\\
       &+\sum\limits_{j=1}^{d_2}\Big[\widetilde{\sigma}_x^{L*j}(t)\widetilde{z}^j(t)
        +\mathbb{E}\big[\widetilde{\sigma}_{\hat{x}}^{L*j}(t)\widetilde{z}^j(t)\big|\mathcal{G}_{1,t}\big]\Big]\\
       &-\sum\limits_{i=1}^n\Big\{\frac{\partial b_x^{L*}}{\partial x_i}(t)q^*(t)p_i(t)\\
       &\qquad\quad+\mathbb{E}\Big[\frac{\partial b_x^{L*}}{\partial\hat{x}_i}(t)q^*(t)p_i(t)\Big|\mathcal{G}_{1,t}\Big]\Big\}\\
       & -\sum\limits_{i=1}^n\Big\{\frac{\partial}{\partial x_i}\bigg(\sum\limits_{j=1}^{d_1}\sigma_x^{L*j}(t)k^j(t)\bigg)p_i(t)\\
       &\quad+\mathbb{E}\Big[\frac{\partial}{\partial\hat{x}_i}\bigg(\sum\limits_{j=1}^{d_1}\sigma_x^{L*j}(t)k^j(t)\bigg)p_i(t)\Big|\mathcal{G}_{1,t}\Big]\Big\}\\
       & -\sum\limits_{i=1}^n\Big\{\frac{\partial}{\partial x_i}\bigg(\sum\limits_{j=1}^{d_1}\widetilde{\sigma}_x^{L*j}(t)\widetilde{k}^j(t)\bigg)p_i(t)\\
       &\quad+\mathbb{E}\Big[\frac{\partial}{\partial\hat{x}_i}
        \bigg(\sum\limits_{j=1}^{d_1}\widetilde{\sigma}_x^{L*j}(t)\widetilde{k}^j(t)\bigg)p_i(t)\Big|\mathcal{G}_{1,t}\Big]\Big\}\\
       &+g_{1xx}^{L*}(t)p(t)+\mathbb{E}\big[g_{1x\hat{x}}^{L*}(t)p(t)\big|\mathcal{G}_{1,t}\big]\\
       &+g_{2x}^{L*}(t)+\mathbb{E}\big[g_{2\hat{x}}^{L*}(t)\big|\mathcal{G}_{1,t}\big]\bigg\}dt\\
       &-z(t)dW(t)-\widetilde{z}(t)d\widetilde{W}(t),\\
  p(0)=&\ 0,\\
  y(T)=&\ G_{1xx}(x^*(T))p(T)+G_{2x}(x^*(T)),
\end{aligned}
\right.
\end{equation}
where we have used $\phi^{L*}(t)\equiv\phi^L\big(t,x^*(t),\hat{x}^*(t),u_2^*(t),\\\hat{u}_2^*(t)\big)$ for $\phi=b,\sigma,\widetilde{\sigma},g_1,g_2$ and all their derivatives.

Now, we have the following two results.
\newpage

\noindent{\bf Proposition 2.3 (Maximum principle of conditional mean-field FBSDE with partial information)}\quad {\it Suppose that {\bf (A1.1)} holds. Let $u_2^*(\cdot)\in\mathcal{U}_2$ be an optimal control for {\bf Problem of the leader} and $(x^*(\cdot),q^*(\cdot),k^*(\cdot),\widetilde{k}^*(\cdot))$ be the optimal state. Let $(y(\cdot),z(\cdot),\widetilde{z}(\cdot),p(\cdot))$ be the adjoint quadruple, then
\begin{equation}\label{maximum condition for the leader-necessary MP}
\begin{aligned}
&\mathbb{E}\bigg[\Big\langle H_{2u_2}\big(t,x^*(t),u_2^*(t),q^*(t),k^*(t),\widetilde{k}^*(t);\\
&\qquad\qquad y(t),z(t),\widetilde{z}(t),p(t)\big),u_2-u_2^*(t)\Big\rangle\\
&\quad+\Big\langle \mathbb{E}\big[H_{2\hat{u}_2}\big(t,x^*(t),u_2^*(t),q^*(t),k^*(t),\widetilde{k}^*(t);y(t),\\
&\qquad\qquad\quad\ z(t),\widetilde{z}(t),p(t)\big)\big|\mathcal{G}_{1,t}\big],\hat{u}_2-\hat{u}_2^*(t)\Big\rangle\bigg|\mathcal{G}_{2,t}\bigg]\\
&\geq0,\ a.e. t\in[0,T], a.s.,
\end{aligned}
\end{equation}
holds for any $u_2\in U_2$.}

{\it Proof.}\quad The maximum condition (\ref{maximum condition for the leader-necessary MP}) can be derived by convex variation and adjoint technique, as Anderson and Djehiche \cite{AD11}. We omit the details for saving space. See also Li \cite{Li12}, Yong \cite{Yong13} and the references therein for mean-field stochastic optimal control problems. \quad$\Box$

\noindent{\bf Proposition 2.4 (Verification theorem of conditional mean-field FBSDE with partial information)}\quad {\it Suppose that {\bf (A1.1)} holds. Let $u_2^*(\cdot)\in\mathcal{U}_2$ and $(x^*(\cdot),q^*(\cdot),k^*(\cdot),\widetilde{k}^*(\cdot))$ be the corresponding state, with $G_{1xx}(x)\equiv G_1\in\mathcal{S}^n$. Let $(y(\cdot),z(\cdot),\widetilde{z}(\cdot),p(\cdot))$ be the adjoint quadruple. For each $(t,\omega)\in[0,T]\times\Omega$, suppose that $H_2\big(t,\cdot,\cdot,\cdot,\cdot,\cdot;y(t),z(t),\widetilde{z}(t),p(t)\big)$ and $G_2(\cdot)$ are convex, and
\begin{equation}\label{maximum condition for the leader-sufficient MP}
\begin{aligned}
 &\mathbb{E}\bigg[H_2\big(t,x^*(t),u_2^*(t),q^*(t),k^*(t),\widetilde{k}^*(t);\\
 &\qquad\qquad y(t),z(t),\widetilde{z}(t),p(t)\big)\\
 &\quad+\mathbb{E}\big[H_2\big(t,x^*(t),u_2^*(t),q^*(t),k^*(t),\widetilde{k}^*(t);y(t),\\
 &\qquad\qquad\quad\ z(t),\widetilde{z}(t),p(t)\big)\big|\mathcal{G}_{1,t}\big]\bigg|\mathcal{G}_{2,t}\bigg]\\
=&\ \max\limits_{u_2\in U_2}\mathbb{E}\bigg[H_2\big(t,x^*(t),u_2,q^*(t),k^*(t),\widetilde{k}^*(t);\\
 &\qquad\qquad\qquad y(t),z(t),\widetilde{z}(t),p(t)\big)\\
 &\quad+\mathbb{E}\big[H_2\big(t,x^*(t),u_2,q^*(t),k^*(t),\widetilde{k}^*(t);y(t),\\
 &\qquad\qquad\quad\ z(t),\widetilde{z}(t),p(t)\big)\big|\mathcal{G}_{1,t}\big]\bigg|\mathcal{G}_{2,t}\bigg],
\end{aligned}
\end{equation}
holds for $a.e. t\in[0,T], a.s.$ Then $u^*_2(\cdot)$ is an optimal control for {\bf Problem of the leader}.}

{\it Proof.}\quad This follows similar to Shi \cite{Shi12}. We omit the details for simplicity. \quad$\Box$

{\bf Remark 2.3}\quad In some applications, for example, the principle-agent problems, the opposite inclusion $\mathcal G_{1,t}\supseteq\mathcal G_{2,t}$ or more complicated relations between them may hold. The study of the leader-follower game in such cases is similar to the one which we studied in this paper. We omit them to limit the length of the paper.

\section{LQ Leader-Follower Stochastic Differential Game with Asymmetric Information}

In order to illustrate the theoretical results in Section 2, we study an LQ leader-follower stochastic differential game with asymmetric information. For notational simplicity, we consider $d_1=d_2=1$. This game is a partially observed model which is a special case of the one in Section 2, but the resulting deduction is technically demanding. We split this section into two subsections, to deal with the problems of the follower and the leader, respectively.

\subsection{The Follower's LQ Problem}

Suppose that the state $(x^{u_1,u_2}(\cdot),\widetilde{x}(\cdot))\in\mathbb{R}^n\times\mathbb{R}^n$ of the system is described by the linear SDE
\begin{equation}\label{state equation for the follower-LQ case}
\left\{
\begin{aligned}
&d\left(\begin{array}{c}x^{u_1,u_2}(t)\\\widetilde{x}(t)\end{array}\right)=\bigg[\left(\begin{array}{cc}A(t)&0\\0&C_1(t)\end{array}\right)\left(\begin{array}{c}x^{u_1,u_2}(t)\\\widetilde{x}(t)\end{array}\right)\\
&\qquad+\left(\begin{array}{cc}B_1(t)&B_2(t)\\0&0\end{array}\right)\left(\begin{array}{c}u_1(t)\\u_2(t)\end{array}\right)\bigg]dt\\
&\qquad+\bigg[\left(\begin{array}{cc}C(t)&0\\0&0\end{array}\right)\left(\begin{array}{c}x^{u_1,u_2}(t)\\\widetilde{x}(t)\end{array}\right)\\
&\qquad+\left(\begin{array}{cc}D_1(t)&D_2(t)\\0&0\end{array}\right)\left(\begin{array}{c}u_1(t)\\u_2(t)\end{array}\right)\\
&\qquad+\left(\begin{array}{c}0\\C_2(t)\end{array}\right)\bigg]dW(t)+\left(\begin{array}{c}\widetilde{C}(t)\\C_3(t)\end{array}\right)d\widetilde{W}(t),\\
&\left(\begin{array}{c}x^{u_1,u_2}(0)\\\widetilde{x}(0)\end{array}\right)=\left(\begin{array}{c}x_0\\\widetilde{x}_0\end{array}\right),
\end{aligned}
\right.
\end{equation}
where $u_1(\cdot)$ takes values in $\mathbb{R}^{m_1},u_2(\cdot)$ takes values in $\mathbb{R}^{m_2}$. Noting that the second part of the state $\widetilde{x}(\cdot)$ is not controlled. Here $A(\cdot),C(\cdot),C_1(\cdot)\in\mathcal{S}^{n},\ \widetilde{C}(\cdot),C_2(\cdot),C_3(\cdot)\in\mathbb{R}^n,\ B_1(\cdot),D_1(\cdot)\in\mathbb{R}^{n\times m_1}$, and $B_2(\cdot),D_2(\cdot)\in\mathbb{R}^{n\times m_2}$ are deterministic, bounded matrix-valued or vector-valued functions, and $x_0,\widetilde{x}_0\in\mathbb{R}^n$.

Suppose that the state $(x^{u_1,u_2}(\cdot),\widetilde{x}(\cdot))$ cannot be directly observed by the follower. Instead, he can observe a related process $Y(\cdot)\in\mathbb{R}$ which is governed by the SDE
\begin{equation}\label{observation equation-LQ case}
\left\{
\begin{aligned}
dY(t)&=\left(\begin{array}{cc}0&h^\top(t)\end{array}\right)\left(\begin{array}{c}x^{u_1,u_2}(t)\\\widetilde{x}(t)\end{array}\right)dt+d\widetilde{W}(t),\\
 Y(0)&=0,
\end{aligned}
\right.
\end{equation}
where $h:[0,T]\rightarrow\mathbb{R}^n$ is a given continuous, bounded vector-valued function.

Inserting (\ref{observation equation-LQ case}) into (\ref{state equation for the follower-LQ case}), we have
\begin{equation}\label{state equation for the follower-LQ case-Y}
\left\{
\begin{aligned}
  dx^{u_1,u_2}(t)&=\big[A(t)x^{u_1,u_2}(t)-h^\top(t)\widetilde{C}(t)\widetilde{x}(t)\\
                 &\quad\ +B_1(t)u_1(t)+B_2(t)u_2(t)\big]dt\\
                 &\quad+\big[C(t)x^{u_1,u_2}(t)+D_1(t)u_1(t)\\
                 &\quad\ +D_2(t)u_2(t)\big]dW(t)+\widetilde{C}(t)dY(t),\\
d\widetilde{x}(t)&=\big[C_1(t)-h^\top(t)C_3(t)\big]\widetilde{x}(t)dt\\
                 &\quad\ +C_2(t)dW(t)+C_3(t)dY(t),\\
   x^{u_1,u_2}(0)&=x_0,\ \widetilde{x}(0)=\widetilde{x}_0.
\end{aligned}
\right.
\end{equation}

Let us now introduce an $\mathbb{R}$-valued process
\begin{equation}\label{P-martingale}
\begin{aligned}
Z(t)&=\exp\left\{-\int_0^th^\top(s)\widetilde{x}(s)d\widetilde{W}(s)\right.\\
    &\qquad\quad\left.-\frac{1}{2}\int_0^t|h^\top(s)\widetilde{x}(s)|^2ds\right\},
\end{aligned}
\end{equation}
which is the solution to the SDE
\begin{equation}\label{R-N derivative}
\left\{
\begin{aligned}
dZ(t)&=-Z(t)h^\top(t)\widetilde{x}(t)d\widetilde{W}(t),\\
 Z(0)&=1.
\end{aligned}
\right.
\end{equation}
It is obvious that $Z(\cdot)$ is a $((\mathcal{F}_t)_{0\leq t\leq T},\mathbb{P})$-martingale, since the Novikov condition holds:
$$
\mathbb{E}\bigg[\exp\Big(\frac{1}{2}\int_0^T\big|h^\top(t)\widetilde{x}(t)\big|^2dt\Big)\bigg]<\infty
$$
(see Example 3, Section 6.2, Liptser and Shiryayev \cite{LS77}). We thus can define a new probability measure $\widetilde{\mathbb{P}}$ such that for any $t,d\widetilde{\mathbb{P}}=Z(t)d\mathbb{P}$, on $\mathcal{F}_t$. By Girsanov's theorem and (\ref{observation equation-LQ case}), $(W(\cdot),Y(\cdot))$ is a 2-dimensional standard Brownian motion defined on $(\Omega,\mathcal{F},(\mathcal{F}_t)_{0\leq t\leq T},\widetilde{\mathbb{P}})$. In this section, we denote $\mathcal{F}_t\equiv\mathcal{F}_t^{W,\widetilde{W}}\equiv\mathcal{F}_t^W\otimes\mathcal{F}_t^{\widetilde{W}}$ and $\mathcal{F}_t^{W,Y}\equiv\mathcal{F}_t^W\otimes\mathcal{F}_t^Y$, where $\mathcal{F}_t^W,\mathcal{F}_t^{\widetilde{W}},\mathcal{F}_t^Y$ are the natural filtration generated by Brownian motions $W(\cdot),\widetilde{W}(\cdot)$ and $Y(\cdot)$, respectively. Moreover, it can be easily checked that $\mathcal{F}_t=\mathcal{F}_t^{W,Y}$.

For chosen $u_2(\cdot)$ by the leader, the follower would like to choose an $\mathcal{F}_t^Y$-adapted control $u_1^*(\cdot)=u_1^*(\cdot;u_2(\cdot))$ to minimize his cost functional
\begin{equation}\label{cost functional for the follower-LQ case}
\begin{aligned}
 &J_1(u_1(\cdot),u_2(\cdot))\\
 &=\frac{1}{2}\mathbb{E}\bigg[\int_0^T\Big(\big\langle Q_1(t)x^{u_1,u_2}(t),x^{u_1,u_2}(t)\big\rangle\\
 &\qquad\qquad\quad+\big\langle N_1(t)u_1(t),u_1(t)\big\rangle\Big)dt\\
 &\qquad\quad+\big\langle G_1x^{u_1,u_2}(T),x^{u_1,u_2}(T)\big\rangle\bigg].
\end{aligned}
\end{equation}
Here $Q_1(\cdot)\in\mathcal{S}^n,N_1(\cdot)\in\mathcal{S}^{m_1}$ are deterministic and bounded matrix-valued functions, and $G_1\in\mathcal{S}^n$. We introduce the following assumption.

\noindent{\bf (A3.1)}\quad{\it $Q_1(t)\geq0,\forall t\in[0,T]$ and $G_1\geq0$.}

The admissible control set $\mathcal{U}_1$ for the follower is defined as (\ref{admissible control set for the follower}) in Section 1, for $\mathcal{G}_{1,t}=\mathcal{F}_t^Y$ in this section. And we will denote the above as {\bf LQ Problem of the follower.}

We apply the maximum principle approach in Subsection 2.1 to solve this partially observed stochastic optimal control problem of the follower. In the following, we will split the process of finding observable optimal control $u_1^*(\cdot)$ for the follower into three steps.

{\it Step 1.} (Optimal control)

First, we put the above problem into the filtered probability space $(\Omega,\mathcal{F},\{\mathcal{F}_t^{W,Y}\}_{0\leq t\leq T},\widetilde{\mathbb{P}})$, in order to apply Propositions 2.1 and 2.2. By (\ref{observation equation-LQ case}) and (\ref{R-N derivative}), we obtain
\begin{equation}\label{equation of Z-by Y}
\left\{
\begin{aligned}
dZ(t)&=Z(t)|h^\top(t)\widetilde{x}(t)|^2dt-Z(t)h^\top(t)\widetilde{x}(t)dY(t),\\
 Z(0)&=1.
\end{aligned}
\right.
\end{equation}
Moreover, we have
\begin{equation}\label{equation of Z-1-by Y}
\left\{
\begin{aligned}
dZ^{-1}(t)&=Z^{-1}(t)h^\top(t)\widetilde{x}(t)dY(t),\\
 Z^{-1}(0)&=1,
\end{aligned}
\right.
\end{equation}
that is,
\begin{equation}\label{Z-1}
\begin{aligned}
Z^{-1}(t)&=\exp\left\{\int_0^th^\top(s)\widetilde{x}(s)dY(s)\right.\\
         &\qquad\quad\left.-\frac{1}{2}\int_0^t|h^\top(s)\widetilde{x}(s)|^2ds\right\}.
\end{aligned}
\end{equation}
By Bayes' formula (Theorem 3.22, Xiong \cite{Xiong08}), we have
\begin{equation}\label{cost functional for the follower-LQ case-by Y}
\begin{aligned}
 &J_1(u_1(\cdot),u_2(\cdot))\\
 &=\frac{1}{2}\widetilde{\mathbb{E}}\bigg[\int_0^T\Big(\big\langle Z^{-1}(t)Q_1(t)x^{u_1,u_2}(t),x^{u_1,u_2}(t)\big\rangle\\
 &\qquad\qquad\quad+\big\langle Z^{-1}(t)N_1(t)u_1(t),u_1(t)\big\rangle\Big)dt\\
 &\qquad\quad+\big\langle Z^{-1}(T)G_1x^{u_1,u_2}(T),x^{u_1,u_2}(T)\big\rangle\bigg],
\end{aligned}
\end{equation}
where $\widetilde{\mathbb{E}}$ is the mathematical expectation with respect
to $\widetilde{\mathbb{P}}$. The Hamiltonian function (\ref{Hamiltonian function for the follower}) of the follower, now takes the form
\begin{equation}\label{Hamiltonian function for the follower-LQ case}
\begin{aligned}
&H_1\big(t,x,\widetilde{x},u_1,u_2;q,k,\widetilde{k};Z^{-1}\big)\\
&=\big\langle q,A(t)x-h^\top(t)\widetilde{C}(t)\widetilde{x}+B_1(t)u_1+B_2(t)u_2\big\rangle\\
&\quad+\big\langle k,C(t)x+D_1(t)u_1+D_2(t)u_2\big\rangle+\big\langle\widetilde{k},\widetilde{C}(t)\big\rangle\\
&\quad-\frac{1}{2}\big\langle Z^{-1}Q_1(t)x,x\big\rangle-\frac{1}{2}\big\langle Z^{-1}N_1(t)u_1,u_1\big\rangle,
\end{aligned}
\end{equation}
and the adjoint equation (\ref{adjoint equation for the follower}) writes
\begin{equation}\label{adjoint equation for the follower-LQ case}
\left\{
\begin{aligned}
-dq(t)=&\ \big[A^\top(t)q(t)+C^\top(t)k(t)\\
       &\quad-Z^{-1}(t)Q_1(t)x^{u_1^*,u_2}(t)\big]dt\\
       &-k(t)dW(t)-\widetilde{k}(t)dY(t),\\
  q(T)=&\ -Z^{-1}(T)G_1x^{u_1^*,u_2}(T),
\end{aligned}
\right.
\end{equation}
with $q(\cdot),k(\cdot),\widetilde{k}(\cdot)$ valued in $\mathbb{R}^n\times\mathbb{R}^n\times\mathbb{R}^n$ being $\mathcal{F}_t^{W,Y}$-adapted processes.

For given $u_2(\cdot)$, suppose that there exists an $\mathcal{F}_t^Y$-adapted optimal control $u_1^*(\cdot)\in\mathcal{U}_1$ for the follower, and the corresponding optimal state $x^{u_1^*,u_2}(\cdot)$ is the solution to (\ref{state equation for the follower-LQ case-Y}), that is,
\begin{equation}\label{state equation-optimal state for the follower}
\left\{
\begin{aligned}
dx^{u_1^*,u_2}(t)&=\big[A(t)x^{u_1^*,u_2}(t)-h^\top(t)\widetilde{C}(t)\widetilde{x}(t)\\
                 &\quad\ +B_1(t)u_1^*(t)+B_2(t)u_2(t)\big]dt\\
                 &\quad+\big[C(t)x^{u_1^*,u_2}(t)+D_1(t)u_1^*(t)\\
                 &\quad\ \ +D_2(t)u_2(t)\big]dW(t)+\widetilde{C}(t)dY(t),\\
 x^{u_1^*,u_2}(0)&=x_0,
\end{aligned}
\right.
\end{equation}
where $\widetilde{x}(\cdot)\in\mathbb{R}^n$ can be solved from the second equation of (\ref{state equation for the follower-LQ case-Y}), since it is a process which is not controlled.

Then (\ref{Hamiltonian function for the follower-LQ case}) together with Proposition 2.1 yields that,
\begin{equation}\label{optimal control for the follower-LQ case}
\begin{aligned}
0=&\ \widehat{Z^{-1}}(t;\widetilde{P})N_1(t)u_1^*(t)-B_1^\top(t)\hat{q}(t;\widetilde{P})\\
  &\ -D_1^\top(t)\hat{k}(t;\widetilde{P}),\ \mbox{for\ } a.e. t\in[0,T],
\end{aligned}
\end{equation}
with
\begin{equation*}
\left\{
\begin{aligned}
&\widehat{Z^{-1}}(t;\widetilde{P}):=\widetilde{\mathbb{E}}[Z^{-1}(t)|\mathcal{F}_t^Y],\\
&\hat{q}(t;\widetilde{P}):=\widetilde{\mathbb{E}}[q(t)|\mathcal{F}_t^Y],\ \hat{k}(t;\widetilde{P}):=\widetilde{\mathbb{E}}[k(t)|\mathcal{F}_t^Y].
\end{aligned}
\right.
\end{equation*}

{\it Step 2.} (Optimal filtering)

Next, noting that the terminal condition of (\ref{adjoint equation for the follower-LQ case}), it is natural to set
\begin{equation}\label{supposed form of q(t)}
q(t)=-Z^{-1}(t)\big[P_1(t)x^{u_1^*,u_2}(t)+\varphi(t)\big],\ t\in[0,T],
\end{equation}
for some deterministic and differentiable $\mathcal{S}^n$-matrix valued function $P_1(\cdot)$, and $\mathbb{R}^n$-valued $\mathcal{F}_t^{W,Y}$-adapted process $\varphi(\cdot)$ admits
\begin{equation}\label{supposed form of varphi(t)}
\left\{
\begin{aligned}
d\varphi(t)=&\ \alpha(t)dt+\beta(t)dY(t),\\
 \varphi(T)=&\ 0.
\end{aligned}
\right.
\end{equation}
In the above, $\alpha(\cdot)\in\mathbb{R}^n$ is $\mathcal{F}_t^{W,Y}$-adapted and $\beta(\cdot)\in\mathbb{R}^n$ is $\mathcal{F}_t^Y$-adapted process, which will be determined later.

Now, applying It\^{o}'s formula to (\ref{supposed form of q(t)}), we have
\begin{equation}\label{applying Ito's formula to q(t)}
\begin{aligned}
dq(t)=&-Z^{-1}(t)\big[P_1(t)dx^{u_1^*,u_2}(t)+\dot{P}_1(t)x^{u_1^*,u_2}(t)dt\\
&\quad+d\varphi(t)\big]-dZ^{-1}(t)\big[P_1(t)x^{u_1^*,u_2}(t)+\varphi(t)\big]\\
&\quad-dZ^{-1}(t)\big[P_1(t)dx^{u_1^*,u_2}(t)+d\varphi(t)\big]\\
=&-Z^{-1}(t)\big[\dot{P}_1(t)x^{u_1^*,u_2}(t)+P_1(t)A(t)x^{u_1^*,u_2}(t)\\
 &\quad+P_1(t)B(t)u_1^*(t)+P_1(t)B_2(t)u_2(t)\\
 &\quad+h^\top(t)\widetilde{x}(t)\beta(t)+\alpha(t)\big]dt\\
 &-Z^{-1}(t)P_1(t)\big[C(t)x^{u_1^*,u_2}(t)\\
 &\quad+D_1(t)u_1^*(t)+D_2(t)u_2(t)\big]dW(t)\\
 &-Z^{-1}(t)\big[P_1(t)\widetilde{C}(t)+P_1(t)h^\top(t)\widetilde{x}(t)x^{u_1^*,u_2}(t)\\
 &\quad+h^\top(t)\widetilde{x}(t)\varphi(t)+\beta(t)\big]dY(t).
\end{aligned}
\end{equation}
Comparing the $\{\cdot\}dW(t)$, $\{\cdot\}dY(t)$ and $\{\cdot\}dt$ terms of (\ref{applying Ito's formula to q(t)}) with those of the adjoint equation (\ref{adjoint equation for the follower-LQ case}), we arrive at
\begin{equation}\label{comparing dW(t)}
\begin{aligned}
k(t)&=-Z^{-1}(t)P_1(t)\big[C(t)x^{u_1^*,u_2}(t)+D_1(t)u_1^*(t)\\
    &\qquad+D_2(t)u_2(t)\big],
\end{aligned}
\end{equation}
\begin{equation}\label{comparing d widetilde W(t)}
\begin{aligned}
\widetilde{k}(t)=&-Z^{-1}(t)\big[P_1(t)\widetilde{C}(t)+P_1(t)h^\top(t)\widetilde{x}(t)x^{u_1^*,u_2}(t)\\
                 &\quad+h^\top(t)\widetilde{x}(t)\varphi(t)+\beta(t)\big],
\end{aligned}
\end{equation}
and
\begin{equation}\label{comparing dt}
\begin{aligned}
\alpha(t)=&\big[-\dot{P}_1(t)-P_1(t)A(t)-A^\top(t)P_1(t)\\
          &\ -C^\top(t)P_1(t)C(t)-Q_1(t)\big]x^{u_1^*,u_2}(t)\\
          &-A^\top(t)\varphi(t)-h^\top(t)\widetilde{x}(t)\beta(t)\\
          &-\big[P_1(t)B_1(t)+C^\top(t)P_1(t)D_1(t)\big]u_1^*(t)\\
          &-\big[P_1(t)B_2(t)+C^\top(t)P_1(t)D_2(t)\big]u_2(t),
\end{aligned}
\end{equation}
respectively. By Bayes' formula, noting (\ref{supposed form of q(t)}), we have
\begin{equation}\label{optimal filtering for q}
\begin{aligned}
\hat{q}(t;\widetilde{P}):=&\ \widetilde{\mathbb{E}}[q(t)|\mathcal{F}_t^Y]=\frac{\mathbb{E}[Z(t)q(t)|\mathcal{F}_t^Y]}{\mathbb{E}[Z(t)|\mathcal{F}_t^Y]}\\
                         =&\ \frac{-P_1(t)\hat{x}^{u_1^*,\hat{u}_2}(t)-\hat{\varphi}(t)}{\mathbb{E}[Z(t)|\mathcal{F}_t^Y]},
\end{aligned}
\end{equation}
where $\hat{x}^{u_1^*,\hat{u}_2}(t):=\mathbb{E}[x^{u_1^*,u_2}(t)|\mathcal{F}_t^Y],\hat{\varphi}(t):=\mathbb{E}[\varphi(t)|\mathcal{F}_t^Y]$. Similarly, by (\ref{comparing dW(t)}) we get
\begin{equation}\label{optimal filtering for k}
\begin{aligned}
&\hat{k}(t;\widetilde{P}):=\widetilde{\mathbb{E}}[k(t)|\mathcal{F}_t^Y]=\frac{\mathbb{E}[Z(t)k(t)|\mathcal{F}_t^Y]}{\mathbb{E}[Z(t)|\mathcal{F}_t^Y]}\\
&\hspace{-5mm}=\frac{-P_1(t)C(t)\hat{x}^{u_1^*,\hat{u}_2}(t)-P_1(t)D_1(t)u_1^*(t)-P_1(t)D_2(t)\hat{u}_2(t)}{\mathbb{E}[Z(t)|\mathcal{F}_t^Y]},
\end{aligned}
\end{equation}
where $\hat{u}_2(t):=\mathbb{E}[u_2(t)|\mathcal{F}_t^Y]$.

Applying Lemma 5.4 in Xiong \cite{Xiong08} to the state equation (\ref{state equation-optimal state for the follower}) of the follower, we derive the filtering equation
\begin{equation}\label{optimal filtering equation for x*u1}
\left\{
\begin{aligned}
d\hat{x}^{u_1^*,\hat{u}_2}(t)&=\big[A(t)\hat{x}^{u_1^*,\hat{u}_2}(t)-h^\top(t)\widetilde{C}(t)\hat{\widetilde{x}}(t)\\
                             &\quad\ +B_1(t)u_1^*(t)+B_2(t)\hat{u}_2(t)\big]dt\\
                             &\quad+\widetilde{C}(t)dY(t),\\
 \hat{x}^{u_1^*,\hat{u}_2}(0)&=\ x_0,
\end{aligned}
\right.
\end{equation}
where $\hat{\widetilde{x}}(t):=\mathbb{E}[\widetilde{x}(t)|\mathcal{F}_t^Y]$ satisfying
\begin{equation}\label{filtering equation tilde x}
\left\{
\begin{aligned}
d\hat{\widetilde{x}}(t)&=\big[C_1(t)-h^\top(t)C_3(t)\big]\hat{\widetilde{x}}(t)dt+C_3(t)dY(t),\\
 \hat{\widetilde{x}}(0)&=\widetilde{x}_0,
\end{aligned}
\right.
\end{equation}
which can be easily solved. Thus (\ref{optimal filtering equation for x*u1}) admits a unique $\mathcal{F}_t^Y$-adapted solution $\hat{x}^{u_1^*,\hat{u}_2}(\cdot)$, for any $\hat{u}_2(\cdot)$.

{\it Step 3.} (Optimal state feedback control)

By (\ref{optimal control for the follower-LQ case}), (\ref{optimal filtering for q}) and (\ref{optimal filtering for k}), in order to obtain state feedback form of the follower's optimal control, the remaining task is to give the representation of $\hat{Z}(t):=\mathbb{E}[Z(t)|\mathcal{F}_t^Y]$, since by Bayes' formula, we have
\begin{equation*}
\begin{aligned}
&\widehat{Z^{-1}}(t;\widetilde{P}):=\widetilde{\mathbb{E}}[Z^{-1}(t)|\mathcal{F}_t^Y]\\
&=\frac{1}{\mathbb{E}[Z(t)|\mathcal{F}_t^Y]}:=\frac{1}{\hat{Z}(t)}\equiv\hat{Z}^{-1}(t).
\end{aligned}
\end{equation*}

Applying Theorem 8.1 of \cite{LS77} to (\ref{R-N derivative}) and (\ref{observation equation-LQ case}), we have
\begin{equation}\label{filtering equation of Z(t)}
\left\{
\begin{aligned}
d\hat{Z}(t)&=-h^\top(t)\hat{Z}(t)\hat{\widetilde{x}}(t)d\hat{\widetilde{W}}(t),\\
 \hat{Z}(0)&=1,
\end{aligned}
\right.
\end{equation}
where $d\hat{\widetilde{W}}(t):=dY(t)-h^\top(t)\hat{\widetilde{x}}(t)dt$.

By (\ref{optimal control for the follower-LQ case}), and supposing that

\noindent{\bf (A3.2)}\quad{\it $\widetilde{N}_1(t)>0$,\ for all $t\in[0,T]$,}

we immediately arrive at
\begin{equation}\label{obsevable optimal control for the follower-LQ case-state feedback form}
\begin{aligned}
&u_1^*(t)\equiv u_1^*\big(t;\hat{x}^{u_1^*,\hat{u}_2}(t),\hat{u}_2(t),\hat{\varphi}(t),\beta(t)\big)\\
&=-\widetilde{N}_1^{-1}(t)\big[\widetilde{S}_1^\top(t)\hat{x}^{u_1^*,\hat{u}_2}(t)+\widetilde{S}(t)\hat{u}_2(t)+B_1^\top(t)\hat{\varphi}(t)\big],
\end{aligned}
\end{equation}
for $a.e.\ t\in[0,T]$, where we denote
\begin{equation*}
\left\{
\begin{aligned}
     \widetilde{N}_1(t):=&\ N_1(t)+D_1^\top(t)P_1(t)D_1(t),\\
     \widetilde{S}_1(t):=&\ P_1(t)B_1(t)+C^\top(t)P_1(t)D_1(t),\\
       \widetilde{S}(t):=&\ D_1^\top(t)P_1(t)D_2(t).
\end{aligned}
\right.
\end{equation*}
Substituting (\ref{obsevable optimal control for the follower-LQ case-state feedback form}) into (\ref{comparing dt}), we can obtain that if the ordinary differential equation
\begin{equation}\label{Riccati equation of P1(t)}
\left\{
\begin{aligned}
 \dot{P}_1(t)&=-P_1(t)A(t)-A^\top(t)P_1(t)-C^\top(t)P_1(t)C(t)\\
             &\quad+\widetilde{S}_1(t)\widetilde{N}_1^{-1}(t)\widetilde{S}_1^\top(t)-Q_1(t),\\
       P_1(T)&=G_1
\end{aligned}
\right.
\end{equation}
admits a unique differentiable solution $P_1(\cdot)$, then
\begin{equation}\label{alpha(t)}
\begin{aligned}
   \alpha(t)&=-\widetilde{S}_1(t)\widetilde{N}_1^{-1}(t)\widetilde{S}_1^\top(t)x^{u_1^*,u_2}(t)\\
            &\quad+\widetilde{S}_1(t)\widetilde{N}_1^{-1}(t)\widetilde{S}_1^\top(t)\hat{x}^{u_1^*,u_2}(t)-A^\top\varphi(t)\\
            &\quad+\big[-A^\top(t)+\widetilde{S}_1(t)\widetilde{N}_1^{-1}(t)B_1^\top(t)\big]\hat{\varphi}(t)\\
            &\quad-h^\top(t)\widetilde{x}(t)\beta(t)-\widetilde{S}_2(t)\big]u_2(t)\\
            &\quad+\big[\widetilde{S}_1(t)\widetilde{N}_1^{-1}(t)\widetilde{S}(t)-\widetilde{S}_2(t)\big]\hat{u}_2(t),
\end{aligned}
\end{equation}
where
\begin{equation*}
             \widetilde{S}_2(t):=P_1(t)B_2(t)+C^\top(t)P_1(t)D_2(t).
\end{equation*}
And the BSDE (\ref{supposed form of varphi(t)}) takes the form
\begin{equation}\label{BSDE of varphi(t)}
\left\{
\begin{aligned}
-d\varphi(t)&=\big[A^\top(t)\varphi(t)+\big(\widetilde{B}_2(t)-B_2(t)\big)^\top\hat{\varphi}(t)\\
            &\quad\ +h^\top(t)\widetilde{x}(t)\beta(t)+\widetilde{\widetilde{S}}_1(t)x^{u_1^*,u_2}(t)\\
            &\quad\ -\widetilde{\widetilde{S}}_1(t)\hat{x}^{u_1^*,\hat{u}_2}(t)+\widetilde{S}_2(t)u_2(t)\\
            &\quad\ +\big(\widetilde{S}_3(t)-\widetilde{S}_2(t)\big)\hat{u}_2(t)\big]dt-\beta(t)dY(t),\\
  \varphi(T)&=0,
\end{aligned}
\right.
\end{equation}
where we set
\begin{equation*}
\left\{
\begin{aligned}
            \widetilde{B}_2(t)&:=B_2(t)-B_1(t)\widetilde{N}_1^{-1}(t)\widetilde{S}(t),\\
\widetilde{\widetilde{S}}_1(t)&:=\widetilde{S}_1(t)\widetilde{N}_1^{-1}(t)\widetilde{S}_1^\top(t),\\
            \widetilde{S}_3(t)&:=-\widetilde{S}_1(t)\widetilde{N}_1^{-1}(t)\widetilde{S}(t)+\widetilde{S}_2(t).
\end{aligned}
\right.
\end{equation*}

We rewrite (\ref{Riccati equation of P1(t)}) as the following
\begin{equation}\label{Riccati equation of P1(t)-full form}
\left\{
\begin{aligned}
 \dot{P}_1(t)&=-P_1(t)A(t)-A^\top(t)P_1(t)\\
             &\quad-C^\top(t)P_1(t)C(t)-Q_1(t)\\
             &\quad+\big[P_1(t)B_1(t)+C^\top(t)P_1(t)D_1(t)\big]\\
             &\qquad\big[N_1(t)+D_1^\top(t)P_1(t)D_1(t)\big]^{-1}\\
             &\qquad\big[P_1(t)B_1(t)+C^\top(t)P_1(t)D_1(t)\big]^\top,\\
       P_1(T)&=G_1,\\
             &\hspace{-9mm}N_1(t)+D_1^\top(t)P_1(t)D_1(t)>0,\ \forall t\in[0,T],
\end{aligned}
\right.
\end{equation}
which is similar to the classical Riccati equation (6.6), Chapter 6 in Yong and Zhou \cite{YZ99}, or (5.2) in \cite{Yong02}. Thus its solvability can be guaranteed.

For given $\hat{u}_2(\cdot)$, plugging the optimal control $u_1^*(\cdot)$ of (\ref{obsevable optimal control for the follower-LQ case-state feedback form}) into (\ref{optimal filtering equation for x*u1}), we have
\begin{equation}\label{optimal filtering equation for x*u1-explicit form}
\left\{
\begin{aligned}
 d\hat{x}^{u_1^*,\hat{u}_2}(t)&=\big[\widetilde{A}(t)\hat{x}^{u_1^*,\hat{u}_2}(t)+\widetilde{S}_4(t)\hat{\varphi}(t)\\
                              &\quad+\widetilde{B}_2(t)\hat{u}_2(t)-h^\top(t)\widetilde{C}(t)\hat{\widetilde{x}}(t)\big]dt\\
                              &\quad+\widetilde{C}(t)dY(t),\\
  \hat{x}^{u_1^*,\hat{u}_2}(0)&=x_0,
\end{aligned}
\right.
\end{equation}
which admits a unique $\mathcal{F}_t^Y$-adapted solution $\hat{x}^{u_1^*,\hat{u}_2}(\cdot)$, where we denote
\begin{equation*}
              \widetilde{S}_4(t):=-B_1(t)\widetilde{N}_1^{-1}(t)B_1^\top(t).
\end{equation*}

Applying Lemma 5.4 in Xiong \cite{Xiong08} to (\ref{BSDE of varphi(t)}), we have
\begin{equation}\label{filtering BSDE of varphi(t)}
\left\{
\begin{aligned}
-d\hat{\varphi}(t)&=\big[\widetilde{A}^\top(t)\hat{\varphi}(t)+h^\top(t)\hat{\widetilde{x}}(t)\beta(t)\\
                  &\quad\ +\widetilde{S}_3(t)\hat{u}_2(t)\big]dt-\beta(t)dY(t),\\
  \hat{\varphi}(T)&=0,
\end{aligned}
\right.
\end{equation}
where
\begin{equation*}
               \widetilde{A}(t):=A(t)-B_1(t)\widetilde{N}_1^{-1}(t)\widetilde{S}_1^\top(t).
\end{equation*}

Putting (\ref{optimal filtering equation for x*u1-explicit form}) and (\ref{filtering BSDE of varphi(t)}) together, we get
\begin{equation}\label{FBSDFE}
\left\{
\begin{aligned}
 d\hat{x}^{u_1^*,\hat{u}_2}(t)&=\big[\widetilde{A}(t)\hat{x}^{u_1^*,\hat{u}_2}(t)+\widetilde{S}_4(t)\hat{\varphi}(t)\\
                              &\quad+\widetilde{B}_2(t)\hat{u}_2(t)-h^\top(t)\widetilde{C}(t)\hat{\widetilde{x}}(t)\big]dt\\
                              &\quad+\widetilde{C}(t)dY(t),\\
            -d\hat{\varphi}(t)&=\big[\widetilde{A}^\top(t)\hat{\varphi}(t)+h^\top(t)\hat{\widetilde{x}}(t)\beta(t)\\
                              &\quad+\widetilde{S}_3(t)\hat{u}_2(t)\big]dt-\beta(t)dY(t),\\
  \hat{x}^{u_1^*,\hat{u}_2}(0)&=x_0,\  \hat{\varphi}(T)=\ 0.
\end{aligned}
\right.
\end{equation}

Note that (\ref{FBSDFE}) is an FBSDFE, whose solvability can be easily obtained, for given $u_2(\cdot)$.

Moreover, we can easily check that the convexity/concavity conditions are satisfied. Then by Proposition 2.2, $u_1^*(\cdot)$ by (\ref{obsevable optimal control for the follower-LQ case-state feedback form}) is really optimal.

We summarize the above in the following theorem.

\noindent{\bf Theorem 3.1}\quad{\it Suppose that assumptions {\bf (A3.1), (A3.2)} hold and the Riccati equation (\ref{Riccati equation of P1(t)-full form}) admits a unique differentiable solution $P_1(\cdot)$. For chosen $u_2(\cdot)$ of the leader, {\bf LQ Problem of the follower} admits an optimal control $u_1^*(\cdot)\equiv u_1^*\big(\cdot;\hat{x}^{u_1^*,\hat{u}_2}(\cdot),\hat{u}_2(\cdot),\hat{\varphi}(\cdot),\beta(\cdot)\big)$ of the state feedback form (\ref{obsevable optimal control for the follower-LQ case-state feedback form}), where process triple $(\hat{x}^{u_1^*,\hat{u}_2}(\cdot),\hat{\varphi}(\cdot),\beta(\cdot))$ is the unique $\mathcal{F}_t^Y$-adapted solution to the FBSDFE (\ref{FBSDFE}).}

\noindent{\bf Remark 3.1}\ Due to the explicit appearance of control $u_2(\cdot)$ in (\ref{adjoint equation for the follower-LQ case}), (\ref{state equation-optimal state for the follower}), we introduce the $\mathcal{F}_t^{W,Y}$-adapted process $\varphi(\cdot)$ in (\ref{supposed form of q(t)}) which satisfies the BSDE (\ref{supposed form of varphi(t)}). This is different from \cite{YZ99} but similar to \cite{Yong02}. More importantly, we will see in next subsection that (\ref{BSDE of varphi(t)}) can be regarded as the new backward ``state" equation of the leader.

{\bf Remark 3.2}\quad Noting that (\ref{FBSDFE}) is an FBSDE system, the value $(\hat{x}^{u_1^*,\hat{u}_2}(t),\hat{\varphi}(t),\beta(t))$ of $(\hat{x}^{u_1^*,\hat{u}_2}(\cdot),\hat{\varphi}(\cdot),\beta(\cdot))$ at time $t$ depends on $\{\hat{u}_2(s);s\in[0,T]\}$. Thus, by (\ref{obsevable optimal control for the follower-LQ case-state feedback form}), $u_1^*(\cdot)$ depends on $\{\hat{u}_2(s);s\in[0,T]\}$. This means that $u_1^*(\cdot)$ is generally anticipating. Thus, it is important to find a ``real" state feedback representation for $u_1^*(\cdot)$ in terms of the original state $x^{u_1^*,u_2}(\cdot)$. This work will be done in the end of this section via the dimensional augmentation introduced by Yong \cite{Yong02}, together with the optimal control $u_2^*(\cdot)$ for the leader.

\subsection{The Leader's LQ Problem}

Now, we are ready to study the stochastic optimal control problem of the leader. Knowing that the follower would take his optimal control $u^*_1(\cdot)\equiv u_1^*\big(\cdot;\hat{x}^{u_1^*,\hat{u}_2}(\cdot),\\\hat{u}_2(\cdot),\hat{\varphi}(\cdot),\beta(\cdot)\big)$ by (\ref{obsevable optimal control for the follower-LQ case-state feedback form}), his state equation (\ref{state equation-optimal state for the follower}) writes
\begin{equation}\label{optimal initial state equation for the leader}
\left\{
\begin{aligned}
 dx^{u_2}(t)&=\Big[A(t)x^{u_2}(t)+\big(\widetilde{A}(t)-A(t)\big)\hat{x}^{\hat{u}_2}(t)\\
            &\quad\ +\widetilde{S}_4(t)\hat{\varphi}(t)+B_2(t)u_2(t)\\
            &\quad\ +\big(\widetilde{B}_2(t)-B_2(t)\big)\hat{u}_2(t)\Big]dt\\
            &\quad+\Big[C(t)x^{u_2}(t)+\widetilde{S}_5(t)\hat{x}^{\hat{u}_2}(t)\\
            &\quad\ +\widetilde{B}^\top_1(t)\hat{\varphi}(t)+D_2(t)u_2(t)\\
            &\quad\ +\widetilde{S}_6(t)\hat{u}_2(t)\Big]dW(t)+\widetilde{C}(t)d\widetilde{W}(t),\\
  x^{u_2}(0)&=x_0,
\end{aligned}
\right.
\end{equation}
where we denote $x^{u_2}(\cdot)\equiv x^{u_1^*,u_2}(\cdot),\hat{x}^{\hat{u}_2}(\cdot)\equiv\hat{x}^{u_1^*,\hat{u}_2}(\cdot)$ for simplicity and
\begin{equation*}
\left\{
\begin{aligned}
            \widetilde{B}_1(t)&:=-B_1(t)\widetilde{N}_1^{-1}(t)D_1^\top(t),\\
            \widetilde{S}_5(t)&:=-D_1(t)\widetilde{N}_1^{-1}(t)\widetilde{S}_1^\top(t),\\
            \widetilde{S}_6(t)&:=-D_1(t)\widetilde{N}_1^{-1}(t)\widetilde{S}(t).
\end{aligned}
\right.
\end{equation*}
The leader would like to choose an $\mathcal{F}_t$-adapted optimal control $u^*_2(\cdot)$ such that his cost functional
\begin{equation}\label{cost functional for the leader-LQ case}
\begin{aligned}
&J_2^L(u_2(\cdot)):=J_2(u_1^*(\cdot),u_2(\cdot))\\
&=\frac{1}{2}\mathbb{E}\bigg[\int_0^T\Big(\big\langle Q_2(t)x^{u_2}(t),x^{u_2}(t)\big\rangle\\
&\qquad\qquad+\big\langle N_2(t)u_2(t),u_2(t)\big\rangle\Big)dt\\
&\qquad\quad+\big\langle G_2x^{u_2}(T),x^{u_2}(T)\big\rangle\bigg]
\end{aligned}
\end{equation}
is minimized, where $Q_2(\cdot)\in\mathcal{S}^n,N_2(\cdot)\in\mathcal{S}^{m_2}$ are deterministic, bounded matrix-valued functions, and $G_2\in\mathcal{S}^n$. We introduce the following assumption.

\noindent{\bf (A3.3)}\quad{\it $Q_2(t)\geq0$,\ for all $t\in[0,T]$ and $G_2\geq0$.}

The admissible control set $\mathcal{U}_2$ of the leader is defined as (\ref{admissible control set for the leader}) in Section 1.2, for $\mathcal{G}_{2,t}=\mathcal{F}_t^{W,\widetilde{W}}$ in this section.

As mentioned in Remark 3.1, the leader has to include the process pair $(\varphi(\cdot),\beta(\cdot))$ as part of his new state processes, since $\hat{\varphi}(\cdot)$ is involved in the coefficients of (\ref{optimal initial state equation for the leader}). Thus for any $u_2(\cdot)$, the state equation of the leader is
\begin{equation}\label{state equation for the leader-LQ case-FBSDE}
\left\{
\begin{aligned}
 dx^{u_2}(t)=&\Big[A(t)x^{u_2}(t)+\big(\widetilde{A}(t)-A(t)\big)\hat{x}^{\hat{u}_2}(t)\\
             &+\widetilde{S}_4(t)\hat{\varphi}(t)+B_2(t)u_2(t)\\
             &+\big(\widetilde{B}_2(t)-B_2(t)\big)\hat{u}_2(t)\Big]dt\\
             &+\Big[C(t)x^{u_2}(t)+\widetilde{S}_5(t)\hat{x}^{\hat{u}_2}(t)\\
             &\quad+\widetilde{B}_1^\top(t)\hat{\varphi}(t)+D_2(t)u_2(t)\\
             &\quad+\widetilde{S}_6(t)\hat{u}_2(t)\Big]dW(t)+\widetilde{C}(t)d\widetilde{W}(t),\\
-d\varphi(t)=&\Big[A^\top(t)\varphi(t)+\big(\widetilde{B}_2(t)-B_2(t)\big)^\top\hat{\varphi}(t)\\
             &+\widetilde{\widetilde{S}}_1(t)x^{u_2}(t)-\widetilde{\widetilde{S}}_1(t)\hat{x}^{\hat{u}_2}(t)\\
             &+\widetilde{S}_2(t)u_2(t)+\big(\widetilde{S}_3(t)-\widetilde{S}_2(t)\big)\hat{u}_2(t)\Big]dt\\
             &-\beta(t)d\widetilde{W}(t),\\
  x^{u_2}(0)=&\ x_0,\ \varphi(T)=0,
\end{aligned}
\right.
\end{equation}
which is a conditional mean-field FBSDE. Its solvability for $\mathcal{F}_t$-adapted solution $(x^{u_2}(\cdot),\varphi(\cdot),\beta(\cdot))$ can be easily guaranteed (since the solvability of FBSDFE (\ref{FBSDFE}) for $(\hat{x}^{u_2}(\cdot),\hat{\varphi}(\cdot),\beta(\cdot))$ has been obtained). And we will denote the above as {\bf LQ Problem of the leader}.

We apply the maximum principle approach in Section 2.2 to solve {\bf LQ Problem of the leader}, which now is a complete information one since $\mathcal{G}_{2,t}\equiv\mathcal{F}_t$. We split this process into three steps.

{\it Step 1.}\ (Optimal control)

Since the process triple $(q(\cdot),k(\cdot),\widetilde{k}(\cdot))$ has been replaced by the new pair $(\varphi(\cdot),\beta(\cdot))$, we define the Hamiltonian function $H_2$ as
\begin{equation}\label{Hamiltonian function for the leader-LQ case}
\begin{aligned}
&\ H_2\big(t,x^{u_2},u_2,\varphi,\beta;y,z,\widetilde{z},p\big)\\
&=\big\langle y,A(t)x^{u_2}+\big(\widetilde{A}(t)-A(t)\big)\hat{x}^{\hat{u}_2}+\widetilde{S}_4(t)\hat{\varphi}\\
&\qquad+B_2(t)u_2+\big(\widetilde{B}_2(t)-B_2(t)\big)\hat{u}_2\big\rangle\\
&\quad+\big\langle z,C(t)x^{u_2}+\widetilde{S}_5(t)\hat{x}^{\hat{u}_2}+\widetilde{B}_1^\top(t)\hat{\varphi}\\
&\qquad\quad+D_2(t)u_2+\widetilde{S}_6(t)\hat{u}_2\big\rangle+\big\langle\widetilde{z},\widetilde{C}(t)\big\rangle\\
&\quad+\big\langle p,A^\top(t)\varphi+\big(\widetilde{B}_2(t)-B_2(t)\big)^\top\hat{\varphi}+\widetilde{\widetilde{S}}_1(t)x^{u_2}\\
&\qquad\quad-\widetilde{\widetilde{S}}_1(t)\hat{x}^{u_2}+\widetilde{S}_2(t)u_2+\big(\widetilde{S}_3(t)-\widetilde{S}_2(t)\big)\hat{u}_2\big\rangle\\
&\quad+\frac{1}{2}\big[\big\langle Q_2(t)x^{u_2},x^{u_2}\big\rangle+\big\langle N_2(t)u_2,u_2\big\rangle\big].
\end{aligned}
\end{equation}
And the adjoint equation (\ref{adjoint equation for the leader}) writes
\begin{equation}\label{adjoint equation for the leader-LQ case}
\left\{
\begin{aligned}
 dp(t)=&\ \Big[A(t)p(t)+\big(\widetilde{B}_2(t)-B_2(t)\big)\hat{p}(t)\\
       &\ +\widetilde{S}_4^\top(t)\hat{y}(t)+\widetilde{B}_1(t)\hat{z}(t)\Big]dt,\\
-dy(t)=&\ \Big[A^\top(t)y(t)+\big(\widetilde{A}(t)-A(t)\big)^\top\hat{y}(t)\\
       &\ +C^\top(t)z(t)+\widetilde{S}_5^\top(t)\hat{z}(t)\\
       &\ +\widetilde{\widetilde{S}}_1(t)p(t)-\widetilde{\widetilde{S}}_1(t)\hat{p}(t)+Q_2(t)x^*(t)\Big]dt\\
       &\ -z(t)dW(t)-\widetilde{z}(t)d\widetilde{W}(t),\\
  p(0)=&\ 0,\ y(T)=G_2x^*(T),
\end{aligned}
\right.
\end{equation}
with
\begin{equation*}
\left\{
\begin{aligned}
                  \hat{p}(t)&:=\mathbb{E}\big[p(t)\big|\mathcal{F}_t^Y\big],\ \hat{y}(t):=\mathbb{E}\big[y(t)\big|\mathcal{F}_t^Y\big],\\
                  \hat{z}(t)&:=\mathbb{E}\big[z(t)\big|\mathcal{F}_t^Y\big],
\end{aligned}
\right.
\end{equation*}
and $(p(\cdot),y(\cdot),z(\cdot))\in\mathbb{R}^n
\times\mathbb{R}^n\times\mathbb{R}^n$ being $\mathcal{F}_t$-adapted processes. Note that (\ref{adjoint equation for the leader-LQ case}) is an FBSDE.

Suppose that there exists an $\mathcal{F}_t$-adapted optimal control $u_2^*(\cdot)\in\mathcal{U}_2$ for the leader, and the corresponding optimal state $(x^*(\cdot),\varphi^*(\cdot),\beta^*(\cdot))\equiv(x^{u_2^*}(\cdot),\varphi^*(\cdot),\beta^*(\cdot))$ is the solution to (\ref{state equation for the leader-LQ case-FBSDE}), that is,
\begin{equation}\label{optimal state equation for the leader-LQ case-FBSDE}
\left\{
\begin{aligned}
 dx^*(t)=&\Big[A(t)x^*(t)+\big(\widetilde{A}(t)-A(t)\big)\hat{x}^*(t)\\
             &+\widetilde{S}_4(t)\hat{\varphi}^*(t)+B_2(t)u_2^*(t)\\
             &+\big(\widetilde{B}_2(t)-B_2(t)\big)\hat{u}_2^*(t)\Big]dt\\
             &+\Big[C(t)x^*(t)+\widetilde{S}_5(t)\hat{x}^*(t)\\
             &\quad+\widetilde{B}_1^\top(t)\hat{\varphi}^*(t)+D_2(t)u_2^*(t)\\
             &\quad+\widetilde{S}_6(t)\hat{u}_2^*(t)\Big]dW(t)+\widetilde{C}(t)d\widetilde{W}(t),\\
-d\varphi^*(t)=&\Big[A^\top(t)\varphi^*(t)+\big(\widetilde{B}_2(t)-B_2(t)\big)^\top\hat{\varphi}^*(t)\\
             &+\widetilde{\widetilde{S}}_1(t)x^*(t)-\widetilde{\widetilde{S}}_1(t)\hat{x}^*(t)\\
             &+\widetilde{S}_2(t)u_2^*(t)+\big(\widetilde{S}_3(t)-\widetilde{S}_2(t)\big)\hat{u}_2^*(t)\Big]dt\\
             &-\beta^*(t)d\widetilde{W}(t),\\
  x^*(0)=&\ x_0,\ \varphi^*(T)=0.
\end{aligned}
\right.
\end{equation}
Then (\ref{Hamiltonian function for the leader-LQ case}) together with Proposition 2.3 yields that
\begin{equation}\label{optimal control for the leader-LQ case}
\begin{aligned}
0=&\ N_2(t)u_2^*(t)+\widetilde{S}_2^\top(t)p(t)+\big(\widetilde{S}_3(t)-\widetilde{S}_2(t)\big)^\top(t)\hat{p}(t)\\
  &\ +B_2^\top(t)y(t)+\big(\widetilde{B}_2(t)-B_2(t)\big)^\top\hat{y}(t)\\
  &\ +D_2^\top(t)z(t)+\widetilde{S}_6^\top(t)\hat{z}(t),\quad a.e.\ t\in[0,T].
\end{aligned}
\end{equation}

Moreover, we can easily check that the convexity condition is satisfied. Then by Proposition 2.4, $u_2^*(\cdot)$ defined by (\ref{optimal control for the leader-LQ case}) is really optimal.

{\it Step 2.}\ (State feedback representation)

The representation of $u_2^*(\cdot)$ through (\ref{optimal control for the leader-LQ case}) is not satisfactory. We expect to obtain its state feedback representation via some Riccati equations. For this target, let us regard the $(x^*(\cdot),p(\cdot))^\top$ as the optimal state and put
\begin{equation}\label{new state for the leader-LQ case}
\begin{aligned}
&X=\left(\begin{array}{c}x^*\\p\end{array}\right),\ \Phi=\left(\begin{array}{c}y\\\varphi^*\end{array}\right),\
Z=\left(\begin{array}{c}z\\0\end{array}\right),\ \widetilde{Z}=\left(\begin{array}{c}\widetilde{z}\\\beta^*\end{array}\right),
\end{aligned}
\end{equation}
and (suppressing some $t$ below)
\begin{equation*}
\left\{
\begin{aligned}
&\mathcal{A}_1:=\left(\begin{array}{cc}A&0\\0&A\end{array}\right),\
 \mathcal{A}_2:=\left(\begin{array}{cc}\widetilde{A}-A&0\\0&\widetilde{B}_2-B_2\end{array}\right),\\
&\mathcal{B}_1:=\left(\begin{array}{cc}0&\widetilde{S}_4\\\widetilde{S}_4^\top&0\end{array}\right),\
 \widetilde{\mathcal{B}}_1:=\left(\begin{array}{cc}0&0\\\widetilde{B}_1&0\end{array}\right),\\
&\mathcal{B}_2:=\left(\begin{array}{c}B_2\\0\end{array}\right),\
 \widetilde{\mathcal{B}}_2:=\left(\begin{array}{c}\widetilde{B}_2-B_2\\0\end{array}\right),\\
&\mathcal{B}_3:=\left(\begin{array}{c}0\\\widetilde{S}_2\end{array}\right),\
 \widetilde{\mathcal{B}}_3:=\left(\begin{array}{c}0\\\widetilde{S}_3-\widetilde{S}_2\end{array}\right),\\
&\mathcal{C}_1:=\left(\begin{array}{cc}C&0\\0&0\end{array}\right),\
 \mathcal{C}_2:=\left(\begin{array}{cc}\widetilde{S}_5&0\\0&0\end{array}\right),\
 \widetilde{\mathcal{C}}:=\left(\begin{array}{c}\widetilde{C}\\0\end{array}\right),\\
&\mathcal{D}_2:=\left(\begin{array}{c}D_2\\0\end{array}\right),\
 \widetilde{\mathcal{D}}_2:=\left(\begin{array}{c}\widetilde{S}_6\\0\end{array}\right),\\
&\mathcal{Q}_2:=\left(\begin{array}{cc}Q_2&\widetilde{\widetilde{S}}_1\\\widetilde{\widetilde{S}}_1&0\end{array}\right),\
 \widetilde{\mathcal{Q}}_2:=\left(\begin{array}{cc}0&-\widetilde{\widetilde{S}}_1\\-\widetilde{\widetilde{S}}_1&0\end{array}\right),\\
&X_0:=\left(\begin{array}{c}x_0\\0\end{array}\right),\
 \mathcal{G}_2:=\left(\begin{array}{cc}G_2&0\\0&0\end{array}\right).
\end{aligned}
\right.
\end{equation*}
Then (\ref{optimal state equation for the leader-LQ case-FBSDE}) with (\ref{adjoint equation for the leader-LQ case}) is equivalent to the following conditional mean-field FBSDE
\begin{equation}\label{optimality system of Problem of the leader-high dimension}
\left\{
\begin{aligned}
    dX(t)=&\ \big[\mathcal{A}_1X(t)+\mathcal{A}_2\hat{X}(t)+\mathcal{B}_1\hat{\Phi}(t)\\
          &\ \ +\widetilde{\mathcal{B}}_1\hat{Z}(t)+\mathcal{B}_2u_2^*(t)+\widetilde{\mathcal{B}}_2\hat{u}_2^*(t)\big]dt\\
          &+\big[\mathcal{C}_1X(t)+\mathcal{C}_2\hat{X}(t)+\widetilde{\mathcal{B}}_1^\top\hat{\Phi}(t)\\
          &\quad+\mathcal{D}_2u_2^*(t)+\widetilde{\mathcal{D}}_2\hat{u}_2^*(t)\big]dW(t)+\widetilde{\mathcal{C}}d\widetilde{W}(t),\\
-d\Phi(t)=&\ \big[\mathcal{Q}_2X(t)+\widetilde{\mathcal{Q}}_2\hat{X}(t)+\mathcal{A}_1^\top\Phi(t)\\
          &\ +\mathcal{A}_2^\top\hat{\Phi}(t)+\mathcal{C}_1^\top Z(t)+\mathcal{C}_2^\top\hat{Z}(t)\\
          &\ +\mathcal{B}_3u_2^*(t)+\widetilde{\mathcal{B}}_3\hat{u}_2^*(t)\big]dt\\
          &-Z(t)dW(t)-\widetilde{Z}(t)d\widetilde{W}(t),\\
     X(0)=&\ X_0,\  \Phi(T)=\mathcal{G}_2X(T),
\end{aligned}
\right.
\end{equation}
where
\begin{equation}\label{filtering estimates-leader-LQ-high dimension}
\left\{
\begin{aligned}
\hat{X}(t)&:=\mathbb{E}\big[X(t)\big|\mathcal{F}_t^Y\big],\ \hat{\Phi}(t):=\mathbb{E}\big[\Phi(t)\big|\mathcal{F}_t^Y\big],\\
\hat{Z}(t)&:=\mathbb{E}\big[Z(t)\big|\mathcal{F}_t^Y\big].
\end{aligned}
\right.
\end{equation}
Supposing that

\noindent{\bf (A3.4)}\quad{\it $N_2^{-1}(t)$ exists,\ for all $t\in[0,T]$,}

then from (\ref{optimal control for the leader-LQ case}), for $a.e.\ t\in[0,T]$, we have
\begin{equation}\label{optimal control for the leader-LQ case-final}
\begin{aligned}
u_2^*(t)&=-N_2^{-1}\big[\mathcal{B}_3^\top X(t)+\widetilde{\mathcal{B}}_3^\top\hat{X}(t)+\mathcal{B}_2^\top\Phi(t)\\
        &\qquad\qquad+\widetilde{\mathcal{B}}_2^\top\hat{\Phi}(t)+\mathcal{D}_2^\top Z(t)+\widetilde{\mathcal{D}}_2^\top\hat{Z}(t)\big],
\end{aligned}
\end{equation}
and
\begin{equation}\label{filtered-optimal control for the leader-LQ case-final}
\begin{aligned}
\hat{u}_2^*(t)&=-N_2^{-1}\big[\big(\mathcal{B}_3+\widetilde{\mathcal{B}}_3\big)^\top\hat{X}(t)+\big(\mathcal{B}_2+\widetilde{\mathcal{B}}_2\big)^\top\hat{\Phi}(t)\\
              &\qquad\qquad+\big(\mathcal{D}_2+\widetilde{\mathcal{D}}_2\big)^\top\hat{Z}(t)\big].
\end{aligned}
\end{equation}
Putting (\ref{optimal control for the leader-LQ case-final}), (\ref{filtered-optimal control for the leader-LQ case-final}) into (\ref{optimality system of Problem of the leader-high dimension}), and letting
\begin{equation*}
\left\{
\begin{aligned}
             \overline{\mathcal{A}}_1&:=\mathcal{A}_1-\mathcal{B}_2N_2^{-1}\mathcal{B}_3^\top,\\
             \overline{\mathcal{A}}_2&:=\mathcal{A}_2-\mathcal{B}_2N_2^{-1}\widetilde{\mathcal{B}}_3^\top-\widetilde{\mathcal{B}}_2N_2^{-1}\big(\mathcal{B}_3+\widetilde{\mathcal{B}}_3\big)^\top,\\
             \overline{\mathcal{B}}_1&:=-\mathcal{B}_2N_2^{-1}\mathcal{B}_2^\top,\quad\overline{\mathcal{B}}_2:=-\mathcal{B}_2N_2^{-1}\mathcal{D}_2^\top,\\
 \overline{\widetilde{\mathcal{B}}}_1&:=\mathcal{B}_1-\mathcal{B}_2N_2^{-1}\widetilde{\mathcal{B}}_2^\top-\widetilde{\mathcal{B}}_2N_2^{-1}\big(\mathcal{B}_2+\widetilde{\mathcal{B}}_2\big)^\top,\\
 \overline{\widetilde{\mathcal{B}}}_2&:=\widetilde{\mathcal{B}}_1-\mathcal{B}_2N_2^{-1}\widetilde{\mathcal{D}}_2^\top-\widetilde{\mathcal{B}}_2N_2^{-1}\big(\mathcal{D}_2+\widetilde{\mathcal{D}}_2\big)^\top,\\
             \overline{\mathcal{C}}_1&:=\mathcal{C}_1-\mathcal{D}_2N_2^{-1}\mathcal{B}_3^\top,\ \overline{\mathcal{D}}_2:=-\mathcal{D}_2N_2^{-1}\mathcal{D}_2^\top,\\
             \overline{\mathcal{C}}_2&:=\mathcal{C}_2-\mathcal{D}_2N_2^{-1}\widetilde{\mathcal{B}}_3^\top-\widetilde{\mathcal{D}}_2N_2^{-1}\big(\mathcal{B}_3+\widetilde{\mathcal{B}}_3\big)^\top,\\
 \overline{\widetilde{\mathcal{D}}}_2&:=-\mathcal{D}_2N_2^{-1}\widetilde{\mathcal{D}}_2^\top-\widetilde{\mathcal{D}}_2N_2^{-1}\big(\mathcal{D}_2+\widetilde{\mathcal{D}}_2\big)^\top,\\
             \overline{\mathcal{Q}}_2&:=\mathcal{Q}_2-\mathcal{B}_3N_2^{-1}\mathcal{B}_3^\top,\\
 \overline{\widetilde{\mathcal{Q}}}_2&:=\widetilde{\mathcal{Q}}_2-\mathcal{B}_3N_2^{-1}\widetilde{\mathcal{B}}_3^\top-\widetilde{\mathcal{B}}_3N_2^{-1}\big(\mathcal{B}_3+\widetilde{\mathcal{B}}_3\big)^\top,
\end{aligned}
\right.
\end{equation*}
we get
\begin{equation}\label{optimality system of Problem of the leader-high dimension-without control}
\left\{
\begin{aligned}
    dX(t)=&\ \big[\overline{\mathcal{A}}_1X(t)+\overline{\mathcal{A}}_2\hat{X}(t)+\overline{\mathcal{B}}_1\Phi(t)\\
          &\quad+\overline{\widetilde{\mathcal{B}}}_1\hat{\Phi}(t)+\overline{\mathcal{B}}_2Z(t)+\overline{\widetilde{\mathcal{B}}}_2\hat{Z}(t)\big]dt\\
          &+\big[\mathcal{C}_1X(t)+\overline{\mathcal{C}}_2\hat{X}(t)+\overline{\mathcal{B}}_2^\top\Phi(t)\\
          &\quad+\overline{\widetilde{\mathcal{B}}}_1^\top\hat{\Phi}(t)+\overline{\mathcal{D}}_2Z(t)+\overline{\widetilde{\mathcal{D}}}_2\hat{Z}(t)\big]dW(t)\\
          &+\widetilde{\mathcal{C}}d\widetilde{W}(t),\\
-d\Phi(t)=&\ \big[\overline{\mathcal{Q}}_2X(t)+\overline{\widetilde{\mathcal{Q}}}_2\hat{X}(t)+\overline{\mathcal{A}}_1^\top\Phi(t)\\
          &\quad+\overline{\mathcal{A}}_2^\top\hat{\Phi}(t)+\mathcal{C}_1^\top Z(t)+\overline{\mathcal{C}}_2^\top\hat{Z}(t)\big]dt\\
          &-Z(t)dW(t)-\widetilde{Z}(t)d\widetilde{W}(t),\\
     X(0)=&\ X_0,\  \Phi(T)=\mathcal{G}_2X(T).
\end{aligned}
\right.
\end{equation}

Up to now, we have every reason to suppose that
\begin{equation}\label{relation of X and Y}
\Phi(t)=\mathcal{P}_1(t)X(t)+\mathcal{P}_2(t)\hat{X}(t)+\mathcal{P}_3(t),
\end{equation}
due to the terminal condition in (\ref{optimality system of Problem of the leader-high dimension-without control}), where $\mathcal{P}_1(\cdot),\mathcal{P}_2(\cdot)$ are both differentiable, deterministic $\mathcal{S}^{2n}$-valued functions with $\mathcal{P}_1(T)=\mathcal{G}_2,\mathcal{P}_2(T)=0$, and $\mathcal{F}_t^{W,Y}$-adapted process $\mathcal{P}_3(\cdot)$ satisfies BSDE
\begin{equation}\label{supposing form of P3(t)}
\left\{
\begin{aligned}
-d\mathcal{P}_3(t)=&\ \lambda(t)dt-\mathcal{Q}_3(t)dY(t),\\
 \mathcal{P}_3(T)=&\ 0.
\end{aligned}
\right.
\end{equation}
In the above, $\lambda(\cdot)\in\mathbb{R}^{2n}$ is $\mathcal{F}_t^{W,Y}$-adapted and $\mathcal{Q}_3(\cdot)\in\mathbb{R}^{2n}$ is $\mathcal{F}_t^Y$-adapted, which will be determined later.

Our remaining task is to decouple the conditional mean-field FBSDE (\ref{optimality system of Problem of the leader-high dimension-without control}), by (\ref{relation of X and Y}) and (\ref{supposing form of P3(t)}). This will lead to a derivation of our system of Riccati equations for $\mathcal{P}_i(\cdot),i=1,2$.

First, by (\ref{optimality system of Problem of the leader-high dimension-without control}) and (\ref{observation equation-LQ case}), we have
\begin{equation}\label{state equation of the leader-high dimension-without control}
\left\{
\begin{aligned}
    dX(t)=&\ \Big[\overline{\mathcal{A}}_1X(t)+\overline{\mathcal{A}}_2\hat{X}(t)+\overline{\mathcal{B}}_1\Phi(t)\\
          &\quad+\overline{\widetilde{\mathcal{B}}}_1\hat{\Phi}(t)+\overline{\mathcal{B}}_2Z(t)+\overline{\widetilde{\mathcal{B}}}_2\hat{Z}(t)\\
          &\quad-\left(\begin{array}{cc}\widetilde{\mathcal{C}}&0\end{array}\right)\left(\begin{array}{c}h^\top(t)\widetilde{x}(t)\\0\end{array}\right)\Big]dt\\
          &+\big[\mathcal{C}_1X(t)+\overline{\mathcal{C}}_2\hat{X}(t)+\overline{\mathcal{B}}_2^\top\Phi(t)\\
          &\quad+\overline{\widetilde{\mathcal{B}}}_1^\top\hat{\Phi}(t)+\overline{\mathcal{D}}_2Z(t)+\overline{\widetilde{\mathcal{D}}}_2\hat{Z}(t)\big]dW(t)\\
          &+\widetilde{\mathcal{C}}dY(t),\\
     X(0)=&\ X_0.
\end{aligned}
\right.
\end{equation}
Noting that
\begin{equation*}
\begin{aligned}
\hat{X}(t;\widetilde{P}):=&\ \widetilde{\mathbb{E}}[\hat{X}(t)|\mathcal{F}_t^Y]=\frac{\mathbb{E}[Z(t)\hat{X}(t)|\mathcal{F}_t^Y]}{\mathbb{E}[Z(t)|\mathcal{F}_t^Y]}\\
                         =&\ \hat{X}(t):=\mathbb{E}[X(t)|\mathcal{F}_t^Y],\ \mbox{etc}.,
\end{aligned}
\end{equation*}
then applying Lemma 5.4 of \cite{Xiong08} to (\ref{state equation of the leader-high dimension-without control}), we get
\begin{equation}\label{filtering state equation of the leader-high dimension-without control}
\left\{
\begin{aligned}
    d\hat{X}(t)=&\ \Big[\big(\overline{\mathcal{A}}_1+\overline{\mathcal{A}}_2\big)\hat{X}(t)+\big(\overline{\mathcal{B}}_1+\overline{\widetilde{\mathcal{B}}}_1\big)\hat{\Phi}(t)\\
                &\ +\big(\overline{\mathcal{B}}_2+\overline{\widetilde{\mathcal{B}}}_2\big)\hat{Z}(t)\\
                &\ -\left(\begin{array}{cc}\widetilde{\mathcal{C}}&0\end{array}\right)\left(\begin{array}{c}h^\top(t)\hat{\widetilde{x}}(t)\\0\end{array}\right)\Big]dt+\widetilde{\mathcal{C}}dY(t),\\
     \hat{X}(0)=&\ X_0.
\end{aligned}
\right.
\end{equation}
Also, by the backward equation of (\ref{optimality system of Problem of the leader-high dimension-without control}), we have
\begin{equation}\label{backward state equation of the leader-high dimension-without control}
\left\{
\begin{aligned}
-d\Phi(t)&=\big[\overline{\mathcal{Q}}_2X(t)+\overline{\widetilde{\mathcal{Q}}}_2\hat{X}(t)+\overline{\mathcal{A}}_1^\top\Phi(t)\\
         &\qquad+\overline{\mathcal{A}}_2^\top\hat{\Phi}(t)+\mathcal{C}_1^\top Z(t)+\overline{\mathcal{C}}_2^\top\hat{Z}(t)\\
         &\qquad+\left(\begin{array}{cc}h^\top(t)\widetilde{x}(t)&0\\0&0\end{array}\right)\widetilde{Z}(t)\big]dt\\
         &\quad-Z(t)dW(t)-\widetilde{Z}(t)dY(t),\\
  \Phi(T)&=\mathcal{G}_2X(T).
\end{aligned}
\right.
\end{equation}
Applying It\^{o}'s formula to (\ref{relation of X and Y}), by (\ref{state equation of the leader-high dimension-without control}), (\ref{filtering state equation of the leader-high dimension-without control}) and (\ref{backward state equation of the leader-high dimension-without control}), we obtain
\begin{equation}\label{Applying Ito's formula to Y}
\begin{aligned}
      &d\Phi(t)=\dot{\mathcal{P}_1}(t)X(t)dt+\mathcal{P}_1(t)dX(t)\\
      &\qquad\quad\ +\dot{\mathcal{P}_2}(t)\hat{X}(t)dt+\mathcal{P}_2(t)d\hat{X}(t)+d\mathcal{P}_3(t)\\
     =&\ \Big\{\big[\dot{\mathcal{P}}_1+\mathcal{P}_1\overline{\mathcal{A}}_1+\mathcal{P}_1\overline{\mathcal{B}}_1\mathcal{P}_1\big]X(t)\\
      &\quad+\big[\dot{\mathcal{P}}_2+\mathcal{P}_2\big(\overline{\mathcal{A}}_1+\overline{\mathcal{A}}_2\big)+\mathcal{P}_2\big(\overline{\mathcal{B}}_1+\overline{\widetilde{\mathcal{B}}}_1\big)\mathcal{P}_1\\
      &\qquad+\mathcal{P}_1\big(\overline{\mathcal{B}}_1+\overline{\widetilde{\mathcal{B}}}_1\big)\mathcal{P}_2+\mathcal{P}_2\big(\overline{\mathcal{B}}_1+\overline{\widetilde{\mathcal{B}}}_1\big)\mathcal{P}_2\\
      &\qquad+\mathcal{P}_1\overline{\mathcal{A}}_2+\mathcal{P}_1\overline{\widetilde{\mathcal{B}}}_1\mathcal{P}_1\big]\hat{X}(t)+\mathcal{P}_1\overline{\mathcal{B}}_2Z(t)\\
      &\quad+\big[\mathcal{P}_1\overline{\widetilde{\mathcal{B}}}_2+\mathcal{P}_2\big(\overline{\mathcal{B}}_2+\overline{\widetilde{\mathcal{B}}}_2\big)\big]\hat{Z}(t)\\
      &\quad+\mathcal{P}_1\overline{\mathcal{B}}_1\mathcal{P}_3(t)-\mathcal{P}_1\left(\begin{array}{cc}\widetilde{\mathcal{C}}&0\end{array}\right)\left(\begin{array}{c}h^\top(t)\widetilde{x}(t)\\0\end{array}\right)\\
      &\quad+\mathcal{P}_1\overline{\widetilde{\mathcal{B}}}_1\hat{\mathcal{P}}_3(t)
       -\mathcal{P}_2\left(\begin{array}{cc}\widetilde{\mathcal{C}}&0\end{array}\right)\left(\begin{array}{c}h^\top(t)\hat{\widetilde{x}}(t)\\0\end{array}\right)\\
      &\quad+\mathcal{P}_2\big(\overline{\mathcal{B}}_1+\overline{\widetilde{\mathcal{B}}}_1\big)\hat{\mathcal{P}}_3(t)-\lambda(t)\Big\}dt\\
      &+\Big\{\big[\mathcal{P}_1\mathcal{C}_1+\mathcal{P}_1\overline{\mathcal{B}}_2^\top\mathcal{P}_1\big]X(t)+\big[\mathcal{P}_1\overline{\mathcal{C}}_2\\
      &\qquad+\mathcal{P}_1\overline{\mathcal{B}}_2^\top\mathcal{P}_2+\mathcal{P}_1\overline{\widetilde{\mathcal{B}}}_1^\top\big(\mathcal{P}_1+\mathcal{P}_2\big)\big]\hat{X}(t)\\
      &\qquad+\mathcal{P}_1\overline{\mathcal{B}}_2^\top\mathcal{P}_3(t)+\mathcal{P}_1\overline{\widetilde{\mathcal{B}}}_1^\top\hat{\mathcal{P}}_3(t)+\mathcal{P}_1\overline{\mathcal{D}}_2Z(t)\\
      &\qquad+\mathcal{P}_1\overline{\widetilde{\mathcal{D}}}_2\hat{Z}(t)\Big\}dW(t)\\
      &\ +\Big\{\big(\mathcal{P}_1+\mathcal{P}_2\big)\widetilde{\mathcal{C}}+\mathcal{Q}_3(t)\Big\}dY(t)\\
     =&-\Big\{\big[\mathcal{Q}_2+\overline{\mathcal{A}}_1^\top\mathcal{P}_1\big]X(t)+\big[\overline{\mathcal{\widetilde{Q}}}_2+\overline{\mathcal{A}}_2^\top\mathcal{P}_1\\
      &\quad\ +\overline{\mathcal{A}}_1^\top\mathcal{P}_2+\overline{\mathcal{A}}_2^\top\mathcal{P}_2\big]\hat{X}(t)+\mathcal{C}_1^\top Z(t)\\
      &\quad\ +\overline{\mathcal{C}}_2^\top\hat{Z}(t)+\left(\begin{array}{cc}h^\top(t)\widetilde{x}(t)&0\\0&0\end{array}\right)\widetilde{Z}(t)\\
      &\quad\ +\overline{\mathcal{A}}_1^\top\mathcal{P}_3(t)+\overline{\mathcal{A}}_2^\top\hat{\mathcal{P}}_3(t)\Big\}dt\\
      &+Z(t)dW(t)+\widetilde{Z}(t)dY(t).
\end{aligned}
\end{equation}
Comparing the coefficients of $\{\cdot\}dW(t)$ and $\{\cdot\}dY(t)$ on both sides of the last $``="$ of (\ref{Applying Ito's formula to Y}), we have
\begin{equation}\label{comparing dW(t)-leader}
\begin{aligned}
  Z(t)&=\big[\mathcal{P}_1\mathcal{C}_1+\mathcal{P}_1\overline{\mathcal{B}}_2^\top\mathcal{P}_1\big]X(t)+\big[\mathcal{P}_1\overline{\mathcal{C}}_2\\
      &\qquad+\mathcal{P}_1\overline{\mathcal{B}}_2^\top\mathcal{P}_2+\mathcal{P}_1\overline{\widetilde{\mathcal{B}}}_1^\top\big(\mathcal{P}_1+\mathcal{P}_2\big)\big]\hat{X}(t)\\
      &\quad+\mathcal{P}_1\overline{\mathcal{D}}_2Z(t)+\mathcal{P}_1\overline{\widetilde{\mathcal{D}}}_2\hat{Z}(t)\\
      &\quad+\mathcal{P}_1\overline{\mathcal{B}}_2^\top\mathcal{P}_3(t)+\mathcal{P}_1\overline{\widetilde{\mathcal{B}}}_1^\top\hat{\mathcal{P}}_3(t),
\end{aligned}
\end{equation}
and
\begin{equation}\label{comparing dY(t)-leader}
\begin{aligned}
  \widetilde{Z}(t)&=\big(\mathcal{P}_1+\mathcal{P}_2\big)\widetilde{\mathcal{C}}+\mathcal{Q}_3(t).
\end{aligned}
\end{equation}
Taking the conditional expectation $\mathbb{E}[\cdot|\mathcal{F}_t^Y]$ on both sides of (\ref{comparing dW(t)-leader}), and supposing that
\begin{equation*}
\begin{aligned}
\hspace{-6mm}{\bf (A3.5)}\ &\widetilde{\mathcal{N}}_2^{-1}:=\big[I_{2n}-\mathcal{P}_1\big(\overline{\mathcal{D}}_2+\overline{\widetilde{\mathcal{D}}}_2(t)\big)\big]^{-1}\\
                 &:=\big[I_{2n}+\mathcal{P}_1\big(\mathcal{D}_2+\widetilde{\mathcal{D}}_2\big)N_2^{-1}\big(\mathcal{D}_2+\widetilde{\mathcal{D}}_2\big)^\top\big]^{-1}\\
                 &\mbox{ exists},\\
\end{aligned}
\end{equation*}
we have
\begin{equation}\label{Z-filtered}
\begin{aligned}
            \hat{Z}(t)=&\ \widetilde{\mathcal{N}}_2^{-1}\Big\{\big[\mathcal{P}_1\big(\mathcal{C}_1+\overline{\mathcal{C}}_2\big)
                       +\mathcal{P}_1\big(\overline{\mathcal{B}}_2+\overline{\widetilde{\mathcal{B}}}_1\big)^\top\\
                       &\times\big(\mathcal{P}_1+\mathcal{P}_2\big)\big]\hat{X}(t)+\mathcal{P}_1\big(\overline{\mathcal{B}}_2+\overline{\widetilde{\mathcal{B}}}_1\big)^\top\hat{\mathcal{P}}_3(t)\Big\}.
\end{aligned}
\end{equation}
Inserting (\ref{Z-filtered}) into (\ref{comparing dW(t)-leader}), and supposing that
\begin{equation*}
\begin{aligned}
\hspace{-1cm}{\bf (A3.6)}\quad\overline{\mathcal{N}}_2^{-1}&:=\big[I_{2n}-\mathcal{P}_1\overline{\mathcal{D}}_2\big]^{-1}\\
&:=\big[I_{2n}+\mathcal{P}_1\mathcal{D}_2N_2^{-1}\widetilde{\mathcal{D}}_2^\top\big]^{-1}\mbox{ exists},
\end{aligned}
\end{equation*}
we have
\begin{equation}\label{Z}
\begin{aligned}
  Z(t)&=\overline{\mathcal{N}}_2^{-1}\Big\{\big[\mathcal{P}_1\mathcal{C}_1+\mathcal{P}_1\overline{\mathcal{B}}_2^\top\mathcal{P}_1\big]X(t)\\
      &\qquad+\Big[\mathcal{P}_1\overline{\mathcal{C}}_2+\mathcal{P}_1\overline{\mathcal{B}}_2^\top\mathcal{P}_2+\mathcal{P}_1\overline{\widetilde{\mathcal{B}}}_1^\top\big(\mathcal{P}_1+\mathcal{P}_2\big)\\
      &\qquad\quad+\mathcal{P}_1\overline{\widetilde{\mathcal{D}}}_2\widetilde{\mathcal{N}}_2^{-1}\big[\mathcal{P}_1\big(\mathcal{C}_1+\overline{\mathcal{C}}_2\big)\\
      &\qquad\quad+\mathcal{P}_1\big(\overline{\mathcal{B}}_2+\overline{\widetilde{\mathcal{B}}}_1\big)^\top\big(\mathcal{P}_1+\mathcal{P}_2\big)\big]\Big]\hat{X}(t)\\
      &\qquad+\mathcal{P}_1\overline{\mathcal{B}}_2^\top\mathcal{P}_3(t)+\big[\mathcal{P}_1\overline{\widetilde{\mathcal{B}}}_1^\top+\mathcal{P}_1\overline{\widetilde{\mathcal{D}}}_2\\
      &\qquad\quad\times\widetilde{\mathcal{N}}_2^{-1}\mathcal{P}_1\big(\overline{\mathcal{B}}_2+\overline{\widetilde{\mathcal{B}}}_1\big)^\top\big]\hat{\mathcal{P}}_3(t)\Big\}.
\end{aligned}
\end{equation}
Comparing the coefficients of $\{\cdot\}dt$ in (\ref{Applying Ito's formula to Y}) and substituting (\ref{comparing dY(t)-leader}), (\ref{Z-filtered}), (\ref{Z}) into it, we obtain the following coupled system of Riccati equations
\begin{equation}\label{system of Riccati equations}
\left\{
\begin{aligned}
      &0=\dot{\mathcal{P}}_1+\mathcal{P}_1\overline{\mathcal{A}}_1+\overline{\mathcal{A}}_1^\top\mathcal{P}_1+\mathcal{P}_1\overline{\mathcal{B}}_1\mathcal{P}_1+\mathcal{Q}_2\\
      &\quad\ +\big(\mathcal{C}_1^\top+\mathcal{P}_1\overline{\mathcal{B}}_2\big)\overline{\mathcal{N}}_2^{-1}\mathcal{P}_1\big(\mathcal{C}_1^\top+\mathcal{P}_1\overline{\mathcal{B}}_2\big)^\top,\\
      &0=\dot{\mathcal{P}}_2+\mathcal{P}_2\big(\mathcal{A}_1+\overline{\mathcal{A}}_2\big)+\big(\mathcal{A}_1+\overline{\mathcal{A}}_2\big)^\top\mathcal{P}_2\\
      &\quad\ +\mathcal{P}_2\big(\overline{\mathcal{B}}_1+\overline{\widetilde{\mathcal{B}}}_1\big)\mathcal{P}_1+\mathcal{P}_1\big(\overline{\mathcal{B}}_1+\overline{\widetilde{\mathcal{B}}}_1\big)\mathcal{P}_2\\
      &\quad\ +\mathcal{P}_2\big(\overline{\mathcal{B}}_1+\overline{\widetilde{\mathcal{B}}}_1\big)\mathcal{P}_2+\mathcal{P}_1\overline{\mathcal{A}}_2+\overline{\mathcal{A}}_2^\top\mathcal{P}_1\\
      &\quad\ +\mathcal{P}_1\overline{\widetilde{\mathcal{B}}}_1\mathcal{P}_1+\overline{\mathcal{\widetilde{Q}}}_2+\big(\mathcal{C}_1^\top+\mathcal{P}_1\overline{\mathcal{B}}_2\big)\\
      &\quad\quad\times\overline{\mathcal{N}}_2^{-1}\mathcal{P}_1\big[\overline{\mathcal{C}}_2+\overline{\widetilde{\mathcal{B}}}_1^\top\mathcal{P}_1
       +\big(\overline{\mathcal{B}}_2+\overline{\widetilde{\mathcal{B}}}_1\big)^\top\mathcal{P}_2\big]\\
      &\quad\ +\big[\big(\mathcal{C}_1^\top+\mathcal{P}_1\overline{\mathcal{B}}_2\big)\overline{\mathcal{N}}_2^{-1}\mathcal{P}_1
       \overline{\widetilde{\mathcal{D}}}_2+\overline{\mathcal{C}}_2^\top+\mathcal{P}_1\overline{\widetilde{\mathcal{B}}}_2\\
      &\qquad\ +\mathcal{P}_2\big(\overline{\mathcal{B}}_2+\overline{\widetilde{\mathcal{B}}}_2\big)\big]\widetilde{\mathcal{N}}_2^{-1}
       \big[\mathcal{P}_1\big(\mathcal{C}_1+\overline{\mathcal{C}}_2\big)\\
      &\qquad\ +\mathcal{P}_1\big(\overline{\mathcal{B}}_2+\overline{\widetilde{\mathcal{B}}}_1\big)^\top\big(\mathcal{P}_1+\mathcal{P}_2\big)\big],\\
      &\mathcal{P}_1(T)=\mathcal{G}_2,\ \mathcal{P}_2(T)=0,
\end{aligned}
\right.
\end{equation}
and
\begin{equation}\label{lambda(t)}
\begin{aligned}
\lambda(t)=&\Big[\overline{\mathcal{A}}_1^\top+\mathcal{P}_1\overline{\mathcal{B}}_1+\big(\mathcal{C}_1^\top+\mathcal{P}_1\overline{\mathcal{B}}_2\big)\overline{\mathcal{N}}_2^{-1}\overline{\mathcal{B}}_2^\top\Big]\mathcal{P}_3(t)\\
           &+\Big[\overline{\mathcal{A}}_2^\top+\mathcal{P}_1\overline{\widetilde{\mathcal{B}}}_1+\mathcal{P}_2\big(\overline{\mathcal{B}}_1+\overline{\widetilde{\mathcal{B}}}_1\big)\\
           &\quad+\big(\mathcal{C}_1^\top+\mathcal{P}_1\overline{\mathcal{B}}_2\big)\overline{\mathcal{N}}_2^{-1}\big[\mathcal{P}_1\overline{\widetilde{\mathcal{B}}}_1^\top+\mathcal{P}_1\overline{\widetilde{\mathcal{D}}}_2\\
           &\qquad\quad\times\widetilde{\mathcal{N}}_2^{-1}\mathcal{P}_1\big(\overline{\mathcal{B}}_2+\overline{\widetilde{\mathcal{B}}}_1\big)^\top\Big]\hat{\mathcal{P}}_3(t)\\
                   &\quad+\left(\begin{array}{cc}h^\top(t)\widetilde{x}(t)&0\\0&0\end{array}\right)\mathcal{Q}_3(t)\\
                   &\quad+\left(\begin{array}{cc}h^\top(t)\big(\widetilde{x}(t)-\hat{\widetilde{x}}(t)\big)&0\\0&0\end{array}\right)\mathcal{P}_2\widetilde{\mathcal{C}}.
\end{aligned}
\end{equation}
Applying Lemma 5.4 in Xiong \cite{Xiong08} to (\ref{supposing form of P3(t)}), we have
\begin{equation}\label{BSDE of filtering P3(t)}
\left\{
\begin{aligned}
-d\hat{\mathcal{P}}_3(t)=&\ \Big\{\Big[\overline{\mathcal{A}}_1^\top+\overline{\mathcal{A}}_2^\top+\big(\mathcal{P}_1+\mathcal{P}_2\big)\big(\overline{\mathcal{B}}_1+\overline{\widetilde{\mathcal{B}}}_1\big)\\
                         &\quad+\big(\mathcal{C}_1^\top+\mathcal{P}_1\overline{\mathcal{B}}_2\big)\overline{\mathcal{N}}_2^{-1}
                          \big[\overline{\mathcal{B}}_2^\top+\mathcal{P}_1\overline{\widetilde{\mathcal{B}}}_1^\top\\
                         &\qquad+\mathcal{P}_1\overline{\widetilde{\mathcal{D}}}_2\widetilde{\mathcal{N}}_2^{-1}\mathcal{P}_1
                          \big(\overline{\mathcal{B}}_2+\overline{\widetilde{\mathcal{B}}}_1\big)^\top\big]\Big]\hat{\mathcal{P}}_3(t)\\
                         &\quad+\left(\begin{array}{cc}h^\top(t)\hat{\widetilde{x}}(t)&0\\0&0\end{array}\right)\mathcal{Q}_3(t)\Big\}dt\\
                         &\quad-\mathcal{Q}_3(t)dY(t),\\
 \hat{\mathcal{P}}_3(T)=&\ 0.
\end{aligned}
\right.
\end{equation}
Noting that the above filtered BSDE (\ref{BSDE of filtering P3(t)}) satisfies the stochastic Lipschitz condition in Wang et al. \cite{WRC07}, then it admits a unique $\mathcal{F}_t^Y$-adapted solution $(\hat{\mathcal{P}}_3(\cdot),\mathcal{Q}_3(\cdot))$. Thus the BSDE (\ref{supposing form of P3(t)}) admits a unique $\mathcal{F}_t^{W,Y}$-adapted solution $(\mathcal{P}_3(\cdot),\mathcal{Q}_3(\cdot))$.

Note that the system of Riccati equations (\ref{system of Riccati equations}) is entirely new. For its solvability, we first rewrite the first equation for $\mathcal{P}_1(\cdot)$ as follows
\begin{equation}\label{system of Riccati equations-P1}
\left\{
\begin{aligned}
      &0=\dot{\mathcal{P}}_1+\mathcal{P}_1\big[\mathcal{A}_1-\mathcal{B}_2N_2^{-1}\mathcal{B}_3^\top\big]\\
      &\quad\ +\big[\mathcal{A}_1-\mathcal{B}_2N_2^{-1}\mathcal{B}_3^\top\big]^\top\mathcal{P}_1-\mathcal{P}_1\mathcal{B}_2N_2^{-1}\mathcal{B}_2^\top\mathcal{P}_1\\
      &\quad\ +\big(\mathcal{C}_1^\top-\mathcal{P}_1\mathcal{B}_2N_2^{-1}\mathcal{D}_2^\top\big)\big[I_{2n}+\mathcal{P}_1\mathcal{D}_2N_2^{-1}\widetilde{\mathcal{D}}_2^\top\big]^{-1}\\
      &\quad\ \times\mathcal{P}_1\big(\mathcal{C}_1^\top-\mathcal{P}_1\mathcal{B}_2N_2^{-1}\mathcal{D}_2^\top\big)^\top+\mathcal{Q}_2,\\
      &\mathcal{P}_1(T)=\mathcal{G}_2.
\end{aligned}
\right.
\end{equation}
Note that $$\big[I_{2n}+\mathcal{P}_1\mathcal{D}_2N_2^{-1}\widetilde{\mathcal{D}}_2^\top\big]^{-1}\mathcal{P}_1=\mathcal{P}_1\big[I_{2n}+\widetilde{\mathcal{D}}_2N_2^{-1}\mathcal{D}_2^\top\mathcal{P}_1\big]^{-1}$$ is symmetric (which can be proved by multiplying both sides by $\big[I_{2n}+\mathcal{P}_1\mathcal{D}_2N_2^{-1}\widetilde{\mathcal{D}}_2^\top\big]$ from left and by $\big[I_{2n}+\widetilde{\mathcal{D}}_2N_2^{-1}\mathcal{D}_2^\top\mathcal{P}_1\big]$ from right). Then by standard Riccati equation theory, (\ref{system of Riccati equations-P1}) admits a unique $\mathcal{S}^{2n}$-valued solution. However, the solvability of the second equation for $\mathcal{P}_2(\cdot)$ is widely open. In this paper, we only impose the solvability of it as an assumption.

Putting (\ref{relation of X and Y}), (\ref{Z-filtered}), (\ref{Z}) into (\ref{optimal control for the leader-LQ case-final}), we obtain that
\begin{equation}\label{state feedback optimal control for the leader-LQ case-final}
\begin{aligned}
&u_2^*(t)=-N_2^{-1}\Big\{\big[\mathcal{B}_3^\top+\mathcal{B}_2^\top\mathcal{P}_1+\mathcal{D}_2^\top\overline{\mathcal{N}}_2^{-1}\mathcal{P}_1\\
&\times\big(\mathcal{C}_1+\overline{\mathcal{B}}_2^\top\mathcal{P}_1\big)\big]X(t)+\Sigma_1\big(\mathcal{P}_1,\mathcal{P}_2\big)\hat{X}(t)\\
&+\Sigma_2\big(\mathcal{P}_1\big)\hat{\mathcal{P}}_3(t)+\big[\mathcal{B}_2^\top+\mathcal{D}_2^\top\overline{\mathcal{N}}_2^{-1}\mathcal{P}_1\overline{\mathcal{B}}_2^\top\big]\mathcal{P}_3(t)\Big\},
\end{aligned}
\end{equation}
for $a.e.\ t\in[0,T]$, where
\begin{equation*}
\left\{
\begin{aligned}
   &\Sigma_1\big(\mathcal{P}_1,\mathcal{P}_2\big):=\widetilde{\mathcal{B}}_3^\top+\mathcal{B}_2^\top\mathcal{P}_2\widetilde{\mathcal{B}}_2^\top\big(\mathcal{P}_1+\mathcal{P}_2\big)\\
         &\quad+\mathcal{D}_2^\top\widetilde{\mathcal{N}}_2^{-1}\Big(\mathcal{P}_1\overline{\mathcal{C}}_2+\mathcal{P}_1\overline{\mathcal{B}}_2^\top\mathcal{P}_2
          +\mathcal{P}_1\overline{\widetilde{\mathcal{B}}}_1^\top\big(\mathcal{P}_1+\mathcal{P}_2\big)\\
         &\quad+\mathcal{P}_1\overline{\widetilde{\mathcal{D}}}_2\widetilde{\mathcal{N}}_2^{-1}\big[\mathcal{P}_1\big(\mathcal{C}_1+\overline{\mathcal{C}}_2\big)\\
         &\quad+\mathcal{P}_1\big(\overline{\mathcal{B}}_2+\overline{\widetilde{\mathcal{B}}}_1\big)^\top\big(\mathcal{P}_1+\mathcal{P}_2\big)\big]\Big)\\
         &\quad+\widetilde{\mathcal{D}}_2^\top\widetilde{\mathcal{N}}_2^{-1}\big[\mathcal{P}_1\big(\mathcal{C}_1+\overline{\mathcal{C}}_2\big)\\
         &\quad+\mathcal{P}_1\big(\overline{\mathcal{B}}_2+\overline{\widetilde{\mathcal{B}}}_1\big)^\top\big(\mathcal{P}_1+\mathcal{P}_2\big)\big],\\
   &\Sigma_2\big(\mathcal{P}_1\big):=\widetilde{\mathcal{B}}_2^\top+\mathcal{D}_2^\top\overline{\mathcal{N}}_2^{-1}\big[\mathcal{P}_1\overline{\widetilde{\mathcal{B}}}_1^\top+\mathcal{P}_1\overline{\widetilde{\mathcal{D}}}_2\\
         &\quad\times\widetilde{\mathcal{N}}_2^{-1}\mathcal{P}_1\big(\overline{\mathcal{B}}_2+\overline{\widetilde{\mathcal{B}}}_1\big)^\top\big]\\
         &\quad+\widetilde{\mathcal{D}}_2^\top\widetilde{\mathcal{N}}_2^{-1}\mathcal{P}_1\big(\overline{\mathcal{B}}_2+\overline{\widetilde{\mathcal{B}}}_1\big)^\top.
\end{aligned}
\right.
\end{equation*}
{\it Step 3.}\ (Optimal state equation)

By (\ref{relation of X and Y}), (\ref{comparing dY(t)-leader}), (\ref{Z-filtered}), (\ref{Z}), we have decoupled the optimality system (\ref{optimality system of Problem of the leader-high dimension-without control}). And the optimal ``state" $X(\cdot)=\big(x^*(\cdot),p(\cdot)\big)^\top$ of the leader should be the $\mathcal{F}_t^{W,Y}$-adapted solution to the conditional mean-field SDE
\begin{equation}\label{close-loop state of the leader}
\left\{
\begin{aligned}
       dX(t)=&\Big\{\Sigma_3(\mathcal{P}_1)X(t)+\Sigma_4(\mathcal{P}_1,\mathcal{P}_2)\hat{X}(t)\\
             &\ +\Sigma_5(\mathcal{P}_1)\mathcal{P}_3(t)+\Sigma_6(\mathcal{P}_1)\hat{\mathcal{P}}_3(t)\Big\}dt\\
             &+\Big\{\Sigma_7(\mathcal{P}_1)X(t)+\Sigma_8(\mathcal{P}_1,\mathcal{P}_2)\hat{X}(t)\\
             &\quad+\Sigma_9(\mathcal{P}_1)\mathcal{P}_3(t)+\Sigma_{10}(\mathcal{P}_1)\hat{\mathcal{P}}_3(t)\Big\}dW(t)\\
             &+\widetilde{\mathcal{C}}d\widetilde{W}(t),\\
       X(0)=&\ X_0,
\end{aligned}
\right.
\end{equation}
where $\hat{X}(\cdot)$ satisfies the filtered SDE
\begin{equation}\label{close-loop state of the follower}
\left\{
\begin{aligned}
       d\hat{X}(t)=&\Big\{\Sigma_{11}(\mathcal{P}_1,\mathcal{P}_2)\hat{X}(t)+\Sigma_{12}(\mathcal{P}_1)\hat{\mathcal{P}}_3(t)\Big\}dt\\
                   &+\widetilde{\mathcal{C}}d\widetilde{W}(t),\\
        \hat{X}(0)=&\ X_0,
\end{aligned}
\right.
\end{equation}
with
\begin{equation*}
\left\{
\begin{aligned}
&\Sigma_3(\mathcal{P}_1):=\overline{\mathcal{A}}_1+\overline{\mathcal{B}}_1\mathcal{P}_1+\overline{\mathcal{B}}_2\overline{\mathcal{N}}_2^{-1}\big[\mathcal{P}_1\mathcal{C}_1+\mathcal{P}_1\overline{\mathcal{B}}_2^\top\mathcal{P}_1\big],\\
&\Sigma_4(\mathcal{P}_1,\mathcal{P}_2):=\overline{\mathcal{A}}_2+\overline{\mathcal{B}}_1\mathcal{P}_2+\overline{\widetilde{\mathcal{B}}}_1\big(\mathcal{P}_1+\mathcal{P}_2\big)\\
             &\quad+\overline{\mathcal{B}}_2\overline{\mathcal{N}}_2^{-1}\Big[\mathcal{P}_1\overline{\mathcal{C}}_2+\mathcal{P}_1\overline{\mathcal{B}}_2^\top\mathcal{P}_2
              +\mathcal{P}_1\overline{\widetilde{\mathcal{B}}}_1^\top\big(\mathcal{P}_1+\mathcal{P}_2\big)\\
             &\quad+\mathcal{P}_1\overline{\widetilde{\mathcal{D}}}_2\widetilde{\mathcal{N}}_2^{-1}\big[\mathcal{P}_1\big(\mathcal{C}_1+\overline{\mathcal{C}}_2\big)
              +\mathcal{P}_1\big(\overline{\mathcal{B}}_2+\overline{\widetilde{\mathcal{B}}}_1\big)^\top\\
             &\quad\times\big(\mathcal{P}_1+\mathcal{P}_2\big)\big]\Big]+\overline{\widetilde{\mathcal{B}}}_2\widetilde{\mathcal{N}}_2^{-1}\big[\mathcal{P}_1\big(\mathcal{C}_1+\overline{\mathcal{C}}_2\big)\\
             &\quad+\mathcal{P}_1\big(\overline{\mathcal{B}}_2+\overline{\widetilde{\mathcal{B}}}_1\big)^\top\big(\mathcal{P}_1+\mathcal{P}_2\big)\big],\\
&\Sigma_5(\mathcal{P}_1):=\overline{\mathcal{B}}_1+\overline{\mathcal{B}}_2\overline{\mathcal{N}}_2^{-1}\mathcal{P}_1\overline{\mathcal{B}}_2^\top,\\
&\Sigma_6(\mathcal{P}_1):=\overline{\widetilde{\mathcal{B}}}_1+\overline{\mathcal{B}}_2\overline{\mathcal{N}}_2^{-1}\big[\mathcal{P}_1\overline{\widetilde{\mathcal{B}}}_1^\top+\mathcal{P}_1\overline{\widetilde{\mathcal{D}}}_2\\
             &\quad\times\widetilde{\mathcal{N}}_2^{-1}\mathcal{P}_1\big(\overline{\mathcal{B}}_2+\overline{\widetilde{\mathcal{B}}}_1\big)^\top\big]
              +\overline{\widetilde{\mathcal{B}}}_2\widetilde{\mathcal{N}}_2^{-1}\mathcal{P}_1\big(\overline{\mathcal{B}}_2+\overline{\widetilde{\mathcal{B}}}_1\big)^\top,\\
&\Sigma_7(\mathcal{P}_1):=\mathcal{C}_1+\overline{\mathcal{B}}_2^\top\mathcal{P}_1+\overline{\mathcal{D}}_2\overline{\mathcal{N}}_2^{-1}\big[\mathcal{P}_1\mathcal{C}_1+\mathcal{P}_1\overline{\mathcal{B}}_2^\top\mathcal{P}_1\big],\\
&\Sigma_8(\mathcal{P}_1,\mathcal{P}_2):=\overline{\mathcal{C}}_2+\overline{\mathcal{B}}_2^\top\mathcal{P}_2+\overline{\widetilde{\mathcal{B}}}_1^\top\big(\mathcal{P}_1+\mathcal{P}_2\big)\\
             &\quad+\overline{\mathcal{D}}_2\overline{\mathcal{N}}_2^{-1}\big[\mathcal{P}_1\overline{\mathcal{C}}_2+\mathcal{P}_1\overline{\mathcal{B}}_2^\top\mathcal{P}_2
              +\mathcal{P}_1\overline{\widetilde{\mathcal{B}}}_1^\top\big(\mathcal{P}_1+\mathcal{P}_2\big)\\
             &\quad+\mathcal{P}_1\overline{\widetilde{\mathcal{D}}}_2\widetilde{\mathcal{N}}_2^{-1}\big[\mathcal{P}_1\big(\mathcal{C}_1+\overline{\mathcal{C}}_2\big)\\
             &\quad+\mathcal{P}_1\big(\overline{\mathcal{B}}_2+\overline{\widetilde{\mathcal{B}}}_1\big)^\top\big(\mathcal{P}_1+\mathcal{P}_2\big)\big]\big]\\
             &\quad+\overline{\widetilde{\mathcal{D}}}_2\widetilde{\mathcal{N}}_2^{-1}\big[\mathcal{P}_1\big(\mathcal{C}_1+\overline{\mathcal{C}}_2\big)\\
             &\quad+\mathcal{P}_1\big(\overline{\mathcal{B}}_2+\overline{\widetilde{\mathcal{B}}}_1\big)^\top\big(\mathcal{P}_1+\mathcal{P}_2\big)\big],\\
&\Sigma_9(\mathcal{P}_1):=\overline{\mathcal{B}}_2^\top+\overline{\mathcal{D}}_2\overline{\mathcal{N}}_2^{-1}\mathcal{P}_1\overline{\mathcal{B}}_2^\top,\\
&\Sigma_{10}(\mathcal{P}_1):=\overline{\widetilde{\mathcal{B}}}_1^\top+\overline{\mathcal{D}}_2\overline{\mathcal{N}}_2^{-1}\big[\mathcal{P}_1\overline{\widetilde{\mathcal{B}}}_1^\top
               +\mathcal{P}_1\overline{\widetilde{\mathcal{D}}}_2\\
             &\quad\times\widetilde{\mathcal{N}}_2^{-1}\mathcal{P}_1\big(\overline{\mathcal{B}}_2+\overline{\widetilde{\mathcal{B}}}_1\big)^\top\big]
             +\overline{\widetilde{\mathcal{D}}}_2\widetilde{\mathcal{N}}_2^{-1}\mathcal{P}_1\big(\overline{\mathcal{B}}_2+\overline{\widetilde{\mathcal{B}}}_1\big)^\top,\\
&\Sigma_{11}(\mathcal{P}_1,\mathcal{P}_2):=\overline{\mathcal{A}}_1+\overline{\mathcal{A}}_2+\big(\overline{\mathcal{B}}_1+\overline{\widetilde{\mathcal{B}}}_1\big)\big(\mathcal{P}_1+\mathcal{P}_2\big)\\
             &\quad+\overline{\mathcal{B}}_2\overline{\mathcal{N}}_2^{-1}\big[\mathcal{P}_1\mathcal{C}_1+\mathcal{P}_1\overline{\mathcal{B}}_2^\top\mathcal{P}_1\big]\\
             &\quad+\overline{\mathcal{B}}_2\overline{\mathcal{N}}_2^{-1}\big[\mathcal{P}_1\overline{\mathcal{C}}_2+\mathcal{P}_1\overline{\mathcal{B}}_2^\top\mathcal{P}_2
              +\mathcal{P}_1\overline{\widetilde{\mathcal{B}}}_1^\top\big(\mathcal{P}_1+\mathcal{P}_2\big)\\
             &\quad+\mathcal{P}_1\overline{\widetilde{\mathcal{D}}}_2\widetilde{\mathcal{N}}_2^{-1}\big[\mathcal{P}_1\big(\mathcal{C}_1+\overline{\mathcal{C}}_2\big)\\
             &\quad+\mathcal{P}_1\big(\overline{\mathcal{B}}_2+\overline{\widetilde{\mathcal{B}}}_1\big)^\top\big(\mathcal{P}_1+\mathcal{P}_2\big)\big]\big]\\
             &\quad+\overline{\widetilde{\mathcal{B}}}_2\widetilde{\mathcal{N}}_2^{-1}\big[\mathcal{P}_1\big(\mathcal{C}_1+\overline{\mathcal{C}}_2\big)\\
             &\quad+\mathcal{P}_1\big(\overline{\mathcal{B}}_2+\overline{\widetilde{\mathcal{B}}}_1\big)^\top\big(\mathcal{P}_1+\mathcal{P}_2\big)\big],\\
&\Sigma_{12}(\mathcal{P}_1):=\overline{\mathcal{B}}_1+\overline{\widetilde{\mathcal{B}}}_1+\overline{\mathcal{B}}_2\overline{\mathcal{N}}_2^{-1}\mathcal{P}_1\overline{\mathcal{B}}_2^\top\\
             &\quad+\overline{\mathcal{B}}_2\overline{\mathcal{N}}_2^{-1}\big[\mathcal{P}_1\overline{\widetilde{\mathcal{B}}}_1^\top+\mathcal{P}_1\overline{\widetilde{\mathcal{D}}}_2\\
             &\quad\times\widetilde{\mathcal{N}}_2^{-1}\mathcal{P}_1\big(\overline{\mathcal{B}}_2+\overline{\widetilde{\mathcal{B}}}_1\big)^\top\big]
              +\overline{\widetilde{\mathcal{B}}}_2\widetilde{\mathcal{N}}_2^{-1}\mathcal{P}_1\big(\overline{\mathcal{B}}_2+\overline{\widetilde{\mathcal{B}}}_1\big)^\top.
\end{aligned}
\right.
\end{equation*}
We summarize the above in the following theorem.

\noindent{\bf Theorem 3.2}\quad{\it Suppose that assumptions {\bf (A3.1)$\sim$\\(A3.6)} hold and the system of Riccati equations (\ref{system of Riccati equations}) admits a differentiable solution pair $(\mathcal{P}_1(\cdot),\mathcal{P}_2(\cdot))$. Let $\hat{X}(\cdot)$ be the $\mathcal{F}_t^Y$-adapted solution to (\ref{close-loop state of the follower}), and $X(\cdot)$ be the $\mathcal{F}_t$-adapted solution to (\ref{close-loop state of the leader}). Define $(Y(\cdot),Z(\cdot),\widetilde{Z}(\cdot))$ by (\ref{relation of X and Y}), (\ref{Z}), (\ref{comparing dY(t)-leader}), respectively. Then (\ref{optimality system of Problem of the leader-high dimension-without control}) holds, where $(\mathcal{P}_3(\cdot),\mathcal{Q}_3(\cdot))$ is the unique $\mathcal{F}_t^{W,Y}$-adapted solution to (\ref{supposing form of P3(t)}). Moreover, the state feedback control $u_2^*(\cdot)$ defined by (\ref{state feedback optimal control for the leader-LQ case-final}) is an optimal control for {\bf Problem of the leader}.}

Noting that the optimal control $u_2^*(\cdot)$ for the leader given by (\ref{state feedback optimal control for the leader-LQ case-final}) is nonanticipating. Likewise, for the follower, the optimal control $u_1^*(t)\equiv u_1^*\big(t;\hat{x}^*(t),\hat{u}_2^*(t),\hat{\varphi}^*(t),\beta^*(t)\big)$ can also be represented in a nonanticipating way. In fact, by (\ref{obsevable optimal control for the follower-LQ case-state feedback form}), noting (\ref{filtered-optimal control for the leader-LQ case-final}), (\ref{relation of X and Y}), (\ref{Z-filtered}) and (\ref{comparing dY(t)-leader}), we have
\begin{equation}\label{state feedback optimal control for the follower-LQ case-final}
\begin{aligned}
&u_1^*(t)=-\widetilde{N}_1^{-1}\big[\widetilde{S}_1^\top\hat{x}^*(t)+\widetilde{S}\hat{u}_2^*(t)\\
         &\qquad\quad+B_1^\top\hat{\varphi}^*(t)+\widetilde{D}_1^\top\beta^*(t)\big]\\
        =&-\widetilde{N}_1^{-1}\Big[\left(\begin{array}{cc}\widetilde{S}_1^\top&0\end{array}\right)\hat{X}(t)+\widetilde{S}\hat{u}_2^*(t)\\
         &\qquad\ +\left(\begin{array}{cc}0&B_1^\top\end{array}\right)\hat{\Phi}(t)
          +\left(\begin{array}{cc}0&\widetilde{D}_1^\top\end{array}\right)\hat{\widetilde{Z}}(t)\Big]\\
        =&-\widetilde{N}_1^{-1}\Big[\left(\begin{array}{cc}\widetilde{S}_1^\top&0\end{array}\right)-\widetilde{S}N_2^{-1}\big[\mathcal{B}_3^\top+\mathcal{B}_2^\top\mathcal{P}_1\\
         &\qquad+\mathcal{D}_2^\top\overline{\mathcal{N}}_2^{-1}\mathcal{P}_1\big(\mathcal{C}_1+\overline{\mathcal{B}}_2^\top\mathcal{P}_1\big)+\Sigma_1\big(\mathcal{P}_1,\mathcal{P}_2\big)\big]\\
         &\qquad+\left(\begin{array}{cc}0&B_1^\top\end{array}\right)\big(\mathcal{P}_1+\mathcal{P}_2\big)\Big]\hat{X}(t)\\
         &-\widetilde{N}_1^{-1}\Big[\left(\begin{array}{cc}0&B_1^\top\end{array}\right)-\widetilde{S}N_2^{-1}\big[\Sigma_2\big(\mathcal{P}_1\big)\\
         &\qquad\qquad+\mathcal{B}_2^\top+\mathcal{D}_2^\top\overline{\mathcal{N}}_2^{-1}\mathcal{P}_1\overline{\mathcal{B}}_2^\top\big]\Big]\hat{\mathcal{P}}_3(t)\\
         &-\widetilde{N}_1^{-1}\left(\begin{array}{cc}0&\widetilde{D}_1^\top\end{array}\right)\big[\big(\mathcal{P}_1+\mathcal{P}_2\big)\widetilde{\mathcal{C}}+\mathcal{Q}_3(t)\big],
\end{aligned}
\end{equation}
for $a.e.\ t\in[0,T]$, which is an observable state feedback representation for the optimal control of the follower, where $\hat{X}(\cdot)$ satisfies (\ref{close-loop state of the follower}) and $(\hat{\mathcal{P}}_3(\cdot),\mathcal{Q}_3(\cdot))$ satisfies (\ref{BSDE of filtering P3(t)}).

Up to now, we have solved our LQ leader-follower stochastic differential game with asymmetric information, and it admits an open-loop Stackelberg equilibrium $(u_1^*(\cdot),u_2^*(\cdot))$. Its state feedback representation is (\ref{state feedback optimal control for the follower-LQ case-final}) and (\ref{state feedback optimal control for the leader-LQ case-final}), respectively. And the corresponding optimal state equation of the leader is (\ref{close-loop state of the leader}), the corresponding optimal state (observable) equation of the follower is (\ref{close-loop state of the follower}).

{\bf Remark 3.4}\quad When we consider the complete information case, i.e., $\widetilde{W}(\cdot)$ disappears and $\mathcal{G}_{1,t}=\mathcal{F}_t$, Theorems 3.1 and 3.2 coincide with Theorems 2.3 and 3.3 in \cite{Yong02}.

\section{Concluding Remarks}

In this paper we have discussed a leader-follower stochastic differential game with asymmetric information, or named a stochastic Stackelberg differential game with asymmetric information. This kind of game problem possesses several attractive features worthy of being highlighted. First, the game problem has the Stackelberg or leader-follower feature, which means the two players act as different roles during the game. Thus the usual approach to deal with game problems, such as Yong \cite{Yong90}, Hamad\`{e}ne \cite{Ha99}, Wu \cite{Wu05}, An and \O ksendal \cite{AO08}, Wang and Yu \cite{WY12}, Yu \cite{Yu12}, Hui and Xiao \cite{HX12,HX14}, Shi \cite{Shi13} where the two players act as equivalent roles, does not apply. Second, the game problem has the asymmetric information between the two players, which was not considered in Yong \cite{Yong02}, \O ksendal et al. \cite{OSU13} and Bensoussan et al. \cite{BCS12}. In detail, the information available to the follower is based on some sub-$\sigma$-algebra of that available to the leader. Stochastic filtering technique is introduced to compute the optimal filtering estimates for the corresponding adjoint processes, which perform as the solution to some FBSDFE. Third, the open-loop Stackelberg equilibrium is represented in its state feedback form for the partial observation case of LQ problem, under some appropriate assumptions on the coefficient matrices in the state equation and the cost functionals. Some new conditional mean-field FBSDEs and system of Riccati equations are first introduced in this paper, to deal with the leader's LQ problem.

Note that in principle, Theorems 3.1 and 3.2 provide a useful tool to seek open-loop Stackelberg equilibrium. As a first step in this direction, we apply our results to LQ models to obtain explicit solutions. We hope to return to the more general case when the states are not linear or the costs are not quadratic in our future research to completely solve the problems in the motivating examples introduced in the introduction. We expect that explicit solutions will not be available for these examples and numerical approximations will be studied. It is worthy to study the closed-loop Stackelberg equilibrium for our problem, as well as the solvability of the system of Riccati equations (\ref{system of Riccati equations}). In addition, many more partially observable cases which are more important and reasonable for applications and technological demanding in its filtering procedure, are highly desirable for further research. These challenging topics will be considered in our future work.

\begin{ack}                               
The authors would like to thank the editor and the anonymous referee for their constructive and insightful comments for improving the quality of
this work. The main content of this paper is presented by the first author in IMS-China International Conference on Statistics and Probability, Kunming, China, June 2015, and International Conference on Mathematical Control Theory-In Memory of Professor Xunjing Li for His 80th Birthday, Chengdu, China, July 2015. Many thanks for discussions and suggestions with Professors Fuzhou Gong, Jin Ma, Jiongmin Yong, Jianfeng Zhang and Yonghui Zhou. The first and second authors would like to thank Department of Mathematics, University of Macau for their hospitality during their visit to Macau.
\end{ack}

\begin{thebibliography}{99}     
\bibitem{AO08}An, T.T.K., \& \O ksendal, B. (2008). Maximum principle for stochastic differential games with partial information. \emph{J. Optim. Theory Appl.}, {\bf 139}(3), 463-483.

\bibitem{AD11}Andersson, D., \& Djehiche, B. (2011). A maximum principle for SDEs of mean-field type. \emph{Appl. Math. Optim.}, {\bf 63}, 341-356.

\bibitem{Basar79}Basar, T.(1979). Stochastic stagewise Stackelberg strategies for linear quadratic systems. In \emph{Stochastic Control Theory and Stochastic Differential Systems}, edited by M. Kohlmann and W. Vogel, Springer-Verlag, Berlin.

\bibitem{BO82}Basar, T., \& Olsder, G.J. (1982). \emph{Dynamic Noncooperative Game Theory}, Academic Press, London.

\bibitem{BCS12}Bensoussan, A., Chen, S., \& Sethi, S.P. (2015). The maximum principle for global solutions of stochastic Stackelberg differential games. \emph{SIAM J. Control Optim.}, {\bf 53}(4), 1956-1981.

\bibitem{BL08}Buckdahn, R., \& Li, J. (2008). Stochastic differential games and viscosity solutions of Hamilton-Bellman-Isaacs equations. \emph{SIAM J. Control Optim.}, {\bf 47}(1), 444-475.

\bibitem{Ha99}Hamad\`{e}ne, S. (1999). Nonzero-sum linear-quadratic stocha-stic differential games and backward-forward equations. \emph{Stoch. Anal. Appl.}, {\bf 17}(1), 117-130.

\bibitem{HPS09}He, X., Prasad, A., \& Sethi, S.P. (2009). Cooperative advertising and pricing in a dynamic stochastic supply chain: feedback Stackelberg strategies. \emph{Prod. Oper. Manag.}, {\bf 18}(1), 78-94.

\bibitem{HWX09}Huang, J., Wang, G., \& Xiong, J. (2009). A maximum principle for partial information backward stochastic control problems with applications, \emph{SIAM J. Control Optim.}, {\bf 48}(4), 2106-2117.

\bibitem{HX12}Hui, Eddie C.M., \& Xiao, H. (2012). Maximum principle for differential games of forward-backward stochastic systems with applications. \emph{J. Math. Anal. Appl.}, {\bf 386}(1), 412-427.

\bibitem{HX14}Hui, Eddie C.M., \& Xiao, H. (2014). Differential games of partial information forward-backward doubly SDE and applications. \emph{ESAIM: COCV}, {\bf 20}(1), 78-94.

\bibitem{Isaacs54}Isaacs, R. (1954-55). \emph{Differential Games}, Parts 1-4. The Rand Corpration, Research Memorandums Nos. RM-1391, RM-1411, RM-1486.

\bibitem{Li12}Li, J. (2012). Stochastic maximum principle in the mean-field controls. \emph{Automatica}, {\bf 48}(2), 366-373.

\bibitem{LS77}Liptser, R.S., \& Shiryayev, A.N. (1977). \emph{Statistics of Random Processes}, Springer-Verlag, New York.

\bibitem{OSU13}\O ksendal, B., Sandal L., \& Ub\o e, J. (2013). Stochastic Stackelberg equilibria with applications to time dependent newsvendor models. \emph{J. Econ. Dyna. \& Control}, {\bf 37}(7), 1284-1299.

\bibitem{Shi12}Shi, J. (2012). Sufficient conditions of optimality for mean-field stochastic control problems. In: \emph{Proc. 12th ICARCV}, Guangzhou, China, 5-7th December, 2012, 747-752.

\bibitem{Shi13}Shi, J. (2013). Relationship between maximum principle and dynamic programming in stochastic differential games and applications. \emph{Amer. J. Oper. Res.}, {\bf 3}(6), 445-453.

\bibitem{Stackelberg34}Stackelberg, H. von. (1934). \emph{Marktform und Gleichgewicht}, Springer, Vienna. (An English translation appeared in \emph{The Theory of the Market Economy}, Oxford University Press, 1952.)

\bibitem{WWX13}Wang, G., Wu, Z., \& Xiong, J. (2013). Maximum principles for forward-backward stochastic control systems with correlated state and observation noises. \emph{SIAM J. Control Optim.}, {\bf 51}(1), 491-524.

\bibitem{WY10}Wang, G. \& Yu, Z. (2010). A Pontryagin's maximum principle for non-zero sum differential games of BSDEs with applications. \emph{IEEE Trans. Autom. Control}, {\bf 55}(7), 1742-1747.

\bibitem{WY12}Wang, G. \& Yu, Z. (2012). A partial information non-zero sum differential game of backward stochastic differential equations with applications. \emph{Automatica}, {\bf 48}(2), 342-352.

\bibitem{WRC07}Wang, J., Ran, Q., \& Chen, Q. (2007). $L^p$ solutions of BSDEs with stochastic Lipschitz condition. \emph{J. Appl. Math. Stoch. Anal.}, {\bf Volume 2007}, Article ID 78196, 14 pages.

\bibitem{Wu05}Wu, Z. (2005). Forward-backward stochastic differential equations, linear quadratic stochastic optimal control and nonzero sum differential games. \emph{J. Syst. Sci. \& Comp.}, {\bf 18}(2), 179-192.

\bibitem{Xiong08}Xiong, J. (2008). \emph{An Introduction to Stochastic Filtering Theory}, Oxford University Press, London.

\bibitem{Yong90}Yong, J. (1990). A zero-sum differential game in a finite duration with switching strategies. \emph{SIAM J. Control Optim.}, {\bf 28}(5), 1234-1250.

\bibitem{Yong02}Yong, J. (2002). A leader-follower stochastic linear quadratic differential games. \emph{SIAM J. Control Optim.}, {\bf 41}(4), 1015-1041.

\bibitem{Yong13}Yong, J. (2013). Linear-quadratic optimal control problems for mean-field stochastic differential equations. \emph{SIAM J.Control Optim.}, {\bf 51}(4), 2809-2838.

\bibitem{YZ99}Yong, J., \& Zhou, X. (1999). \emph{Stochastic Controls: Hamiltonian Systems and HJB Equations}, Springer-Verlag, New York.

\bibitem{Yu12}Yu, Z. (2012). Linear-quadratic optimal control and nonzero-sum differential game of forward-backward stochastic system. \emph{Asian J. Control}, {\bf 14}(1), 173-185.

\bibitem{Zhang90}Zhang, Q. (1990). Controlled partially observed diffusions with correlated noise. \emph{Appl. Math. Optim.}, {\bf 22}(3), 265-285.


\end{thebibliography}

\end{document}